\font\gothic=eufm10.
.
\font\sets=msbm10.
\font\stampatello=cmcsc10.
\font\symbols=msam10.
\def\0{{\bf 0}}
\def\1{{\bf 1}}
\def\defineq{\buildrel{def}\over{=}}
\def\definiz{\buildrel{def}\over{\Longleftrightarrow}}
\def\C{\hbox{\sets C}}
\def\H{\hbox{\sets H}}
\def\N{\hbox{\sets N}}

\def\Primes{\hbox{\sets P}}
\def\Z{\hbox{\sets Z}}
\def\square{\hbox{\vrule\vbox{\hrule\phantom{s}\hrule}\vrule}}
\def\supporto{{\rm supp}}
\def\EssBdd{\hbox{\symbols n}\,}
\def\WA{(\hbox{\stampatello WA})}
\def\WSA{(\hbox{\stampatello WSA})}
\def\WWA{(\hbox{\stampatello WWA})}
\def\DH{(\hbox{\stampatello DH})}
\def\DD{(\hbox{\stampatello DD})}
\def\BH{(\hbox{\stampatello BH})}
\def\HL{(\hbox{\stampatello H-L})}
\def\REEF{(\hbox{\stampatello R.E.E.F.})}
\def\Reef{\hbox{\stampatello REEF}}
\def\FAI{(\hbox{\stampatello FAI})}
\def\KVL{(\hbox{\stampatello KVL})}
\def\KVI{(\hbox{\stampatello KVI})}
\def\WIN{(\hbox{\stampatello WIN})}

\def\Irr{{\hbox{\rm Irr}}}
\def\IrrPdF{{\hbox{\rm Irr}^{(P)}_d\,F}}
\def\IrrP1F{{\hbox{\rm Irr}^{(P)}_1\,F}}
\def\IrrPF{{\hbox{\rm Irr}^{(P)}\,F}}

\def\IrrQFdF{{\hbox{\rm Irr}^{(Q_F)}_d\,F}}
\def\IrrQPdF{{\hbox{\rm Irr}^{(P_F)}_d\,F}}
\def\SingSer{\hbox{\gothic S}}
\def\Carmichael{{\rm Car}}
\def\Wintner{{\rm Win}}
\def\CarmichaelT{\Carmichael \; }
\def\WintnerT{\Wintner \; }
\def\CarmichaelTP{\Carmichael^{(P)} \; }
\def\WintnerTP{\Wintner^{(P)} \; }
\def\CarmichaelPqF{\Carmichael^{(P)}_q F\;}
\def\WintnerPqF{\Wintner^{(P)}_q F\;}
\def\CloudF{<F>}
\def\SmoothCloudF{\subset F\supset}
\def\SmoothCloud0{\subset {\bf 0}\supset}
\def\AbsCloudF{\CloudF_{\hbox{\stampatello abs}}}
\def\Rvl{(\hbox{\stampatello Rvl})}
\def\NSL{(\hbox{\stampatello NSL})}
\def\IPP{(\hbox{\stampatello IPP})}
\def\ETD{(\hbox{\stampatello ETD})}

\def\qed{\hfill $\square$\par}
\def\nondivide{\!\!\not\,\mid }
\def\sumflat{\mathop{{\sum}^{\flat}}}
\par
\centerline{\bf A smooth summation of Ramanujan expansions}
\bigskip
\centerline{Giovanni Coppola}\footnote{ }{MSC $2010$: $11{\rm N}05$, $11{\rm P}32$, $11{\rm N}37$ - Keywords: Ramanujan expansion, correlation, $2k-$twin primes} 
\bigskip
\bigskip
\bigskip
\par
\noindent
{\bf Abstract}. 
We studied {\it Ramanujan series} \enspace $\sum_{q=1}^{\infty}G(q)c_q(a)$, where $c_q(a)$ is the well-known Ramanujan sum and the complex numbers $G(q)$, as $q\in \N$, are the {\it Ramanujan coefficients}; of course, we mean, implicitly, that the series converges pointwise, in all natural $a$, as its partial sums $\sum_{q\le Q}G(q)c_q(a)$ converge in $\C$, when $Q\to \infty$. Motivated by our recent study of infinite and finite {\it Euler products} for the Ramanujan series, in which we assumed $G$ multiplicative, we look at a kind of (partial) smooth summations. These are \enspace $\sum_{q\in (P)}G(q)c_q(a)$, where the indices $q$ in $(P)$ means that all prime factors $p$ of $q$ are up to $P$ (fixed); then, we pass to the limit over $P\to \infty$. Notice that this kind of partial sums over $P-${\it smooth numbers} (i.e., in $(P)$, see the above) make up an infinite sum, themselves, $\forall P\in \Primes$ \thinspace fixed, in general; however, our summands contain $c_q(a)$, that has a {\it vertical limit}, i.e. it's supported over indices $q\in \N$ for which the $p-$adic valuations of, resp., $q$ and $a$, namely $v_p(q)$, resp., $v_p(a)$ satisfy \enspace $v_p(q)\le v_p(a)+1$ \enspace and this is true $\forall p\le P$ ($P$'s fixed). 
\par
\noindent
In other words, $\forall G:\N \rightarrow \C$, here, \enspace $\sum_{q\in (P)}G(q)c_q(a)$ \enspace is a {\stampatello finite sum}, $\forall a\in \N$, $\forall P\in \Primes$ fixed: we will call \enspace $\sum_{q=1}^{\infty}G(q)c_q(a)$ \enspace a {\it Ramanujan smooth series} if and only if \enspace $\exists \lim_P \sum_{q\in (P)}G(q)c_q(a)\in \C$, $\forall a\in \N$. 
\par
\noindent
Notice a {\it very important property} : Ramanujan smooth series and Ramanujan series need not to be the same. 
\par
\noindent
We prove : Ramanujan smooth series converge under {\it Wintner Assumption}. (This is not necessarily true for Ramanujan series.) We apply this to correlations and to the Hardy--Littlewood \lq \lq $2k$-Twin Primes\rq \rq Conjecture. 

\bigskip
\bigskip

\par
\noindent{\bf 1. Introduction. Main results for: arithmetic functions, correlations and $2k-$twin primes} 
\bigskip
\par
\noindent
We pursue our study of Ramanujan expansions with smooth moduli, started in [C1]. There, we obtained pointwise converging Ramanujan expansions, for some arithmetic functions having Eratosthenes transform supported over smooth numbers : say, the $F:\N \rightarrow \C$ with \lq \lq smooth divisors\rq \rq; we then applied this general result (see [C1], Theorem 1), to the correlations satisfying a reasonable hypothesis (see [C1], Corollary 1). 
\par
Here, a new kind of summing Ramanujan expansions will give us a \lq \lq new world\rq \rq, of elementary results about convergence and, notably, for more general arithmetic functions (no restriction on their divisors, here). 
\par 
Among these, following Theorem 1, a completely unexpected, new version, say, of Delange Theorem [De] about the convergence of Ramanujan expansions : if we confine to the summation of partial sums on smooth numbers, we can get their convergence, but with a weaker hypothesis with respect to Delange's (i.e., $\DH$, following) and this, actually, is the assumption in the Wintner's Criterion (i.e., $(2.1)$ in $(ii)$, see Theorem 2.1 in Chapter VIII of [ScSp]); that we'll call the {\it Wintner assumption}, abbreviated $\WA$, for $F :\N \rightarrow \C$ having [W] {\it Eratosthenes transform} $F'\defineq F\ast \mu$, where $\mu$ is M\"obius function and $\ast$ is Dirichlet product [T]: 
\vskip-0.10in
$$
\sum_{d=1}^{\infty}{{|F'(d)|}\over d}<\infty. 
\leqno{\WA}
$$
\vskip-0.10in
\par
\noindent
However, $\WA$ is not sufficient for the convergence of classical partial sums (see $\S5.1$). In fact, for this we need Delange Hypothesis (compare $(6)$ in [De]), next $\DH$; in which $\omega(d)\defineq |\{ p\in\Primes, p|d\}|$ is the number of prime factors of $d\in\N$ (whence, see [T], $2^{\omega(d)}=\sum_{t|d}\mu^2(t)$ is the number of {\it square-free} divisors of $d$) : 
\vskip-0.10in
$$
\sum_{d=1}^{\infty}{{2^{\omega(d)}|F'(d)|}\over d}<\infty. 
\leqno{\DH}
$$
\vskip-0.10in
\par
\noindent
Like we did in [C1], we write \enspace $(V)\defineq \{n\in \N : (n,p)=1, \forall p>V\}$ \enspace for the set of $V-${\it smooth numbers}, while \enspace $)V( \defineq \{n\in \N : (n,p)=1, \forall p\le V\}$ \enspace is the set of $V-${\it sifted numbers}. Notice that : $(V)\cap \;)V(\;=\{1\}$, $\forall V\in\N$. 
\par
\noindent
We write $V=P\in \Primes$ hereafter, so that $(P)$ and $)P($ avoid the trivial case $(1)=)1(=\{1\}$. 
\bigskip
\par				
\noindent
In the following, we use the classical notation\enspace $\ll$\enspace of Vinogradov ($A\ll B$ means $|A|\le C \cdot B$, for some constant $C>0$), with\enspace $\ll_{\varepsilon}$\enspace indicating a dependence on\enspace $\varepsilon>0$, arbitrarily small usually, in the\enspace $\ll$\enspace constant. 
\par
As usual, we say that $F:\N \rightarrow \C$ satisfies the {\it Ramanujan Conjecture}, by definition, when : $\forall \varepsilon>0$, $\exists C=C(\varepsilon)$ : $|F(n)|\le C\cdot n^{\varepsilon}$, $\forall n\in \N$ (large enough), i.e., in Vinogradov notation, 
$$
\forall \varepsilon>0, 
\qquad
F(n)\ll_{\varepsilon} n^{\varepsilon},
\quad
\hbox{\rm as } \enspace n\to \infty.
\leqno{(\hbox{\stampatello Ramanujan\enspace Conjecture})}
$$
\par
\noindent
(In other papers, we write $F\EssBdd 1$ for that, also calling $F$ \lq \lq {\stampatello essentially bounded}\rq \rq: compare [C3] and [CM].)
\par
We rely, here and in [C1], on the fact that all $F:\N \rightarrow \C$ satisfying Ramanujan Conjecture and having $F'$ supported on smooth numbers, say $(P)$, have a nice behavior for the convergence issues related to Ramanujan expansions and their coefficients. This is based, at last, on the following bound (compare [C1], Lemma 3, for all the details), in which\enspace $\varepsilon>0$\enspace is arbitrarily small : 
$$
\sum_{m\in (P)} m^{\varepsilon-1}=\prod_{p\le P}\sum_{K=0}^{\infty}\left(p^{\varepsilon-1}\right)^K
=\prod_{p\le P}{1\over {1-p^{\varepsilon-1}}}
<\infty, 
\leqno{(1)}
$$
\par
\noindent
and notice that the same series, but without the condition \lq \lq $m\in (P)$\rq \rq, of course, is a diverging one. 
\par
This elementary estimate (coming from multiplicativity of $m^{\varepsilon-1}$, w.r.t. $m\in\N$) seems to be not so powerful; however, it implies that $F$ satisfying Ramanujan Conjecture, with $P-$smooth divisors, satisfy Delange Hypothesis (see $\DH$ above), that (thanks to [De] main result) implies : Carmichael coefficients $\CarmichaelT F$, see the following, equal Wintner coefficients $\WintnerT F$, see the following (compare Theorem 1 in [C1]).
\par
The $\WA$ is called Wintner Assumption, because Wintner [W] was the first to work with it for the Ramanujan expansions; first of all, by positivity it implies the existence of all the \lq \lq {\it Wintner coefficients}\rq \rq, say, of our $F$, namely : 
$$
\Wintner_q F\defineq \sum_{d\equiv 0\bmod q}{{F'(d)}\over d}, 
\quad
\forall q\in \N,
$$
\par
\noindent
converging (even absolutely) from $\WA$; this also implies the existence of all the following limits in all the, say, \lq \lq {\it Carmichael coefficients}\rq \rq, of our $F$, where $c_q(n)$ is the {\it Ramanujan sum} [R] of modulus $q$ \& argument $n$ (we recall soon after), namely: 
$$
\Carmichael_q F\defineq {1\over {\varphi(q)}}\lim_x {1\over x}\sum_{n\le x}F(n)c_q(n), 
\quad
\forall q\in \N,
$$
\par
\noindent
with $\varphi(q)\defineq |\{ n\le q : (n,q)=1\}|$ the {\it Euler totient function}. Wintner [W] proved that $\WA \Rightarrow \CarmichaelT F=\WintnerT F$, namely $\Carmichael_q F=\Wintner_q F$, $\forall q\in\N$ here. If all these $q-$coefficients exist, these two, $\CarmichaelT F$, resp., $\WintnerT F$, may be called, resp., {\stampatello Carmichael Transform}, resp., {\stampatello Wintner Transform}, of our $F$. (Of course, existence implies uniqueness, for both these transforms; that are arithmetic functions, themselves.)
\par
The {\stampatello Ramanujan smooth expansion} of our $F$, where $c_q(a)\defineq \sum_{j\le q,(j,q)=1}\cos{{2\pi ja}\over q}$ is the well-known {\it Ramanujan sum} [R],[M], holds with these coefficients, under $\WA$ (see next Theorem 1): $\forall a\in \N$, fixed, 
$$
F(a)=\lim_P \sum_{q\in (P)}(\Carmichael_q F)c_q(a)
=\lim_P \sum_{q\in (P)}(\Wintner_q F)c_q(a). 
$$
\medskip
\par
\noindent
We will call hereafter \enspace ${\displaystyle \sum_{q=1}^{\infty}G(q)c_q(a) }$ \enspace a {\stampatello Ramanujan smooth series}, say, of coefficient $G:\N \rightarrow \C$, by definition, when the limit \enspace ${\displaystyle \lim_P \sum_{q\in (P)}G(q)c_q(a) }$ exists in $\C$, for all natural $a$. 
\par
A big warning is that the classical Ramanujan series, defined if \enspace ${\displaystyle \exists \lim_Q \sum_{q\le Q}G(q)c_q(a)\in\C }$, is A PRIORI different from this. (Compare $\S5.1$ for the example of $G$ {\it constant}.) 
\par
Thus, all the results we consider (Lemmas, Theorems \& Corollaries) are about this \lq \lq smooth summation\rq \rq. 
\par
\noindent
We write: \thinspace $\supporto(F)\defineq \{n\in \N : F(n)\neq 0\}$ \thinspace the {\it support} of any \thinspace $F:\N \rightarrow \C$. \enspace We start with our main results. 

\vfill
\eject

\par				
\noindent{\bf 1.1. General Theorems for arithmetic functions}
\bigskip
We give a kind of improvement, of Delange main result [De], inasmuch our partial sums are smooth : instead of $\DH$, we need $\WA$. \enspace In the following, \enspace QED\enspace is the end of a part of a Proof, ending with a \hfill $\square$
\smallskip
\par
\noindent {\bf Theorem 1}. ({\stampatello Wintner's \lq \lq Dream Theorem\rq \rq}) {\it Let } $F:\N\rightarrow\C$ {\it satisfy Wintner Assumption $\WA$. Then} 
$$
\forall a\in \N, 
\quad 
F(a)=\lim_P \sum_{q\in (P)}\left(\Wintner_q F\right)c_q(a)
 =\lim_P \sum_{q\in (P)}\left(\Carmichael_q F\right)c_q(a). 
$$ 
\par
\noindent
{\it The additional hypothesis that } $\WintnerT F$ {\it is, say, smooth supported: } $\supporto(\WintnerT F)\subseteq (Q)$ {\it for some prime } $Q$, {\it gives}
$$
\forall a\in \N, 
\quad 
F(a)=\sum_{q\in (Q)}\left(\Wintner_q F\right)c_q(a)
 =\sum_{q\in (Q)}\left(\Carmichael_q F\right)c_q(a). 
$$ 
\par
\noindent
{\it In particular, in case } $\supporto(\WintnerT F)$ {\it is finite, say}
$$
\exists Q\in \N \enspace : \enspace 
\Wintner_q F =0, 
\qquad
\forall q>Q,
$$
{\it we have} 
$$
\forall a\in \N, 
\quad 
F(a)=\sum_{q\le Q}\left(\Wintner_q F\right)c_q(a)
 =\sum_{q\le Q}\left(\Carmichael_q F\right)c_q(a). 
$$ 
\par
\noindent {\bf Proof}. Fix $a\in \N$, take $P\ge a$, $P\in \Primes$, getting from Lemma 1, $(3)$, 
$$
F(a)=\sum_{d\in (P)}{{F'(d)}\over d}\sum_{{q\in (P)}\atop {q|d}}c_q(a); 
$$
\par
\noindent
then Lemma 1, $(4)$, together with Wintner assumption gives the following double series absolute convergence:  
$$
\sum_{d\in (P)}{{|F'(d)|}\over d}\sum_{q|d}|c_q(a)|\le \sum_{d=1}^{\infty}{{|F'(d)|}\over d}\sum_{q\in (P)}|c_q(a)|
<\infty, 
$$
\par
\noindent
allowing the exchange of these $d,q$ sums: 
$$
F(a)=\sum_{q\in (P)}\Big(\sum_{{d\in (P)}\atop {d\equiv 0\bmod q}}{{F'(d)}\over d}\Big)c_q(a) 
 =\sum_{q\in (P)}\left(\Wintner_q F\right)c_q(a) - \sum_{q\in (P)}\Big(\sum_{{d\not \in (P)}\atop {d\equiv 0\bmod q}}{{F'(d)}\over d}\Big)c_q(a); 
$$
\par
\noindent
another exchange, for these other sums, is possible for the same reason: 
$$
\sum_{q\in (P)}\Big(\sum_{{d\not \in (P)}\atop {d\equiv 0\bmod q}}{{F'(d)}\over d}\Big)c_q(a)
=\sum_{d\not \in (P)}{{F'(d)}\over d}\sum_{{q\in (P)}\atop {q|d}}c_q(a), 
$$
\par
\noindent
implying, from Lemma 1, $(5)$: 
$$
\Big|\sum_{q\in (P)}\Big(\sum_{{d\not \in (P)}\atop {d\equiv 0\bmod q}}{{F'(d)}\over d}\Big)c_q(a)\Big|
\le 
\sum_{d\not \in (P)}{{|F'(d)|}\over d}\Big|\sum_{{q\in (P)}\atop {q|d}}c_q(a)\Big|
 \le a\sum_{d>P}{{\left|F'(d)\right|}\over d}
 \buildrel{P}\over{\longrightarrow} 0, 
$$
\par
\noindent
completing the first part.\enspace QED\enspace The case of $\supporto(\WintnerT F)\subseteq (Q)$ follows from : $(Q)\subseteq (P)$, $\forall P>Q$.\hfill QED 
\par
\noindent
In particular, when $\supporto(\WintnerT F)\subseteq [1,Q]$, use : $q\le Q$ $\Rightarrow$ $q\in (Q)$ and previous case.\qed 
\medskip
\par
\noindent {\bf Remark 1.} A shorter alternative proof (see $\S3$) follows from Lemma $2$.\hfill $\diamond$ 
\medskip
\par				
In what follows, we use the expression {\it fixed length Ramanujan expansion} to indicate a finite Ramanujan expansion \enspace $\sum_{q\le Q}G(q)c_q(a)$ \enspace where $Q\in \N$ is an absolute constant (not $a-$dependent, in particular). In this paper, we will not use the expression {\it finite Ramanujan expansion}, since its length may depend on $a\in \N$.  
\medskip
\par
\noindent {\bf Remark 2.} We proved in [C3] : $F$ has a fixed length Ramanujan expansion $\Longleftrightarrow $ $\supporto(F')$ is finite.\hfill $\diamond$ 
\medskip
\par
In [CM], we gave many characterizations, for correlations $F(a)=C_{f,g}(N,a)\defineq \sum_{n\le N}f(n)g(n+a)$, of the condition: $\supporto(F')$ is finite. 
\medskip
\par
Even if we are considering smooth partial sums, in case they have fixed length, of course, they are the same of classical partial sums. In other words, fixed length partial sums, of course, do converge in any of the summation methods we choose! This trivial remark is applied, in next result : it characterizes the finiteness of Ramanujan series partial sums, say, whenever the Wintner coefficients are, in turn, finitely supported.  

\bigskip

\par
\noindent
We recall the notation $\0(n)\defineq 0$, $\forall n\in \N$, for the {\it null-function}.  
\medskip
\par
Two new characterizations arise, for FIXED LENGTH Ramanujan expansions; first one is a little bit technical, in next result, where second equivalence implies: $F'$ finitely supported $\Leftrightarrow$ $F'$ \lq \lq {\it smooth-supported}\rq \rq. 
\smallskip
\par
\noindent {\bf Theorem 2}. {\it Let } $F:\N\rightarrow\C$ {\it have finite } $\supporto(\WintnerT F)$. {\it Then}
$$
\supporto(F')\enspace \hbox{\it is\enspace finite}
\qquad
\Longleftrightarrow
\qquad
\lim_P \sum_{{r\in )P(}\atop {r>1}}{{F'(dr)}\over r}=\0(d)
\qquad
\Longleftrightarrow
\qquad
\exists Q\in \N \enspace : \enspace \supporto(F')\subseteq (Q).
$$
\par
\noindent {\bf Remark 3.} The series here is defined as 
$$
\sum_{{r\in )P(}\atop {r>1}}{{F'(dr)}\over r}\defineq \lim_x \sum_{{r\in )P(}\atop {1<r\le x}}{{F'(dr)}\over r}, 
$$
\par
\noindent
that exists in $\C$ when $\exists \WintnerT F$, as proved in Lemma 2, $\S2$ (there, compare Remark 5).\hfill $\diamond$ 
\smallskip
\par
\noindent {\bf Proof}. We prove the first and the second equivalence in both directions, considering a large prime $P$. 
\par
\noindent
Since $r\in )P($ and $r>1$ implies $r>P$, whence $dr>P$, $\forall d\in \N$, first \lq \lq $\Rightarrow$\rq \rq \thinspace follows.\hfill QED
\par
\noindent
From $(7)$ of Lemma 2, also first \lq \lq $\Leftarrow$\rq \rq \thinspace follows.\hfill QED 
\par
\noindent
Second \lq \lq $\Leftarrow$\rq \rq \thinspace follows from: $r\in )P($ and $r>1$ $\Rightarrow$ $\exists p>P, p|r$ $\Rightarrow$ $F'(dr)=0$, $\forall d\in \N$.\hfill QED
\par
\noindent
Finally, second \lq \lq $\Rightarrow$\rq \rq \thinspace follows from first \lq \lq $\Leftarrow$\rq \rq \thinspace and the triviality: $\supporto(F')\subseteq [1,Q]$ $\Rightarrow$ $\supporto(F')\subseteq (Q)$.\qed

\bigskip

\par
We give an important \lq \lq summary\rq \rq, for sufficient conditions to get Ramanujan smooth expansions, with Wintner coefficients.  
\smallskip
\par
\noindent {\bf Remark 4.} Let $F:\N\rightarrow\C$ have $\WintnerT F$. If at least one of the following three hypotheses holds:
$$
F'\enspace \hbox{\rm has\enspace finite\enspace support}
\qquad \hbox{\rm OR} \qquad 
F'\enspace \hbox{\rm has\enspace smooth\enspace support}
\qquad \hbox{\rm OR} \qquad 
F\enspace \hbox{\rm satisfies\enspace Wintner\enspace Assumption}, 
$$
\par
\noindent
then $F$ has a Ramanujan smooth expansion, with Wintner coefficients. We give an immediate justification, for this. 
From Theorem 1, $\WA$ $\Rightarrow $ the thesis, while finite support implies, trivially, smooth support, too; then we restrict to smoothness of $F'$ support: use Lemma 3, $\S3$ and Remark 6.\hfill $\diamond$
\medskip
\par
Actually, Lemma 3 in $\S3$ gives an equivalent condition for the Ramanujan smooth expansion, with Wintner coefficients. 
\medskip
\par
\noindent
In the forthcoming subsections we present:
\smallskip
\par
\item{$\diamond$} in next subsection, an application to \lq \lq correlations\rq \rq, that satisfy a \lq \lq reasonable hypothesis\rq \rq;
\smallskip
\par
\item{$\diamond$} then, in subsection $1.3$, a particular, but noteworthy case of \lq \lq reasonable correlation\rq \rq: the $2k-$twin primes correlation, in Hardy-Littlewood Conjecture; this is proved under Wintner Assumption (giving a new Conditional Proof stronger than the one we gave in  [C0], under Delange Hypothesis). 

\vfill
\eject

\par				
\noindent
A short glance to the following sections:
\medskip
\par
\item{$\diamondsuit$} Section 2, \lq \lq Lemmata for the Theorems\rq \rq, supplies the Lemmas for Theorems 1 \& 2 Proofs: Lemma 1 gives elementary calculations; while, Lemma 2 is the core of present paper: it presents a kind of \lq \lq arithmetic orthogonality\rq \rq, realizing {\stampatello Wintner's $P-$orthogonality Decomposition}, after a decomposition into two orthogonal sets of indices, namely, the $P-$smooth and the $P-$sifted (here $P$ is any fixed prime). 
\medskip
\par
\item{$\diamondsuit$} Section 3, \lq \lq A deeper look into Ramanujan smooth expansions: Ramanujan-Wintner smooth expansions\rq \rq, gives a characterization of arithmetic functions having the Ramanujan smooth expansion, with Wintner coefficients, in Lemma 3. Also, it provides a shorter Proof for Theorem 1. 
\medskip
\par
\item{$\diamondsuit$} Section 4, starting from an idea in [C1], gives \lq \lq {\stampatello local expansions}\rq \rq \enspace which have {\stampatello P}$-$smooth coefficients (both Wintner's \& Carmichael's) that converge to the coefficients in {\stampatello Ramanujan} smooth {\stampatello expansions}, compare \lq \lq Theorem 1(Smooth Version)\rq \rq. (A kind of stronger Theorem 1, under $\WSA$, Wintner's Smooth Assumption, weaker than Wintner Assumption.) The properties of these $P-$smooth coefficients are then studied in three sets of Arithmetic Functions.
\hfill
\break
From 8th-version onwards, we add Properties 1 and 2. 
\medskip
\par
\item{$\diamondsuit$} Section 5, continuing to expose \& generalize our elementary methods. Speaking about: \lq \lq {\stampatello Ramanujan Clouds}\rq \rq; generalizations of Wintner Assumption (like the $\WSA$, in $\S4$, quoted above, and beyond) \& of the $\REEF$, that we introduce for correlations in next $\S1.2$; and further generalizations: of the $\Reef$ for arithmetic functions $F$ with finite support for $\WintnerT \,F$, a kind of decomposition for $F$ in two parts that are {\stampatello analytic} (an entire function!) and {\stampatello irregular} (from \lq \lq Irregular Series\rq \rq), a brief study of irregular series of multiplicative functions. 
\hfill
\break
Then, we deepen two important issues, expanding previous version5. 
\hfill
\break
First, $\S5.6$, we study the Counterexample 1 in third version of [C1], that proves: $\BH$ for correlations doesn't imply the $\REEF$, providing some interesting details for this very simple correlation.
\hfill
\break
We add, from version8 onwards, Curiosity 1.
\hfill
\break
Second, $\S5.7$, we explicitly calculate $P-$smooth Carmichael-Wintner coefficients for the imaginary exponentials, whence for $\BH-$correlations, proving that they all converge to classical Carmichael-Wintner coefficients, as $P\to \infty$ in primes. A very important difference, from version 6 to 7, is a correction, i.e. $q''$ definition. 
\medskip
\par
\item{$\diamondsuit$} Last but not least: a glance at Euler products, links between Eratosthenes \& Wintner Transforms \lq \lq in Wintner's style\rq \rq, [W], with further Remarks, and a brief coming soon for future work, are in Section 6. Version 9 adds new results, \lq \lq Crossing {\stampatello Horizontal and Vertical Limits}\rq \rq, in 6.3, to get the $\Reef$.  

\bigskip
\bigskip

\par
\noindent{\bf 1.2. Applications for the correlations satisfying Basic Hypothesis}
\bigskip
\par
\noindent
Given two arithmetic functions $f,g:\N \rightarrow \C$, for their {\stampatello correlation}\enspace ${\displaystyle C_{f,g}(N,a)\defineq \sum_{n\le N}f(n)g(n+a) }$, that has Eratosthenes Transform\enspace ${\displaystyle C_{f,g}'(N,t)\defineq \sum_{a|t}C_{f,g}(N,a)\mu(t/a) }$, we assume [C1] the {\stampatello Basic Hypothesis}: 
$$
g(m)\defineq \sum_{q|m,q\le Q}g'(q), 
\enspace \forall m\in \N, \enspace
\hbox{\it with} \enspace Q\le N,
\quad \hbox{\it and} \quad
C_{f,g}(N,a)\enspace \hbox{\it is}\enspace \hbox{\stampatello fair},
\leqno{\BH}
$$
\par
\noindent
where the condition to be {\stampatello fair} for \enspace $C_{f,g}(N,a)$ \enspace means that the dependence on $a$ is only in the argument of $g(n+a)$ (not inside $f$, nor in $g$). The main consequences are given in Proposition 1 of [C1]; in particular, $\BH$ for $C_{f,g}(N,a)$ implies that $C_{f,g}'(N,d)$ satisfies Ramanujan Conjecture (from the boundedness of $C_{f,g}(N,a)$) and Carmichael-Wintner coefficients (i.e., Carmichael \& Wintner coefficients are the same) of $C_{f,g}(N,a)$ are 
$$
{{\widehat{g}(q)}\over {\varphi(q)}}\sum_{n\le N}f(n)c_q(n), 
\enspace \hbox{\rm with} \enspace 
\widehat{g}(q)\defineq \sum_{d\le Q, d\equiv 0\bmod q}{{g'(d)}\over d},
\quad
\forall q\in\N. 
$$
\par				
\noindent
(Since $\;\supporto(\widehat{g})\subseteq [1,Q]$, also for these coefficients $[1,Q]$ contains their support: outside $[1,Q]$ they vanish !) 
\par
\noindent
The Ramanujan expansion with these coefficients (given in $(iii)$ of Theorem 1 [CM]) is called (see [C3], $\S4$) the {\stampatello Ramanujan exact explicit formula} : 
$$
C_{f,g}(N,a)=\sum_{q\le Q}\left({{\widehat{g}(q)}\over {\varphi(q)}}\sum_{n\le N}f(n)c_q(n)\right)c_q(a), 
\quad 
\forall a\in\N. 
\leqno{\REEF}
$$
\par
We start with our first application, the strongest, for the Correlations. 
\smallskip
\par
\noindent {\bf Corollary 1}. ({\stampatello the $\REEF$ follows from Basic Hypothesis and Wintner Assumption})
\par
\noindent
{\it Let the correlation $C_{f,g}(N,a)$ satisfy } $\BH$ {\it and } $\WA$. {\it Then, the $\REEF$ holds.} 
\smallskip
\par
\noindent {\bf Proof}. From $(iii)$ of Proposition 1 in [C1], $\BH$ gives the finitely-supported Carmichael-Wintner coefficients above. Apply Theorem 1 to $F(a)=C_{f,g}(N,a)$.\qed

\bigskip

New characterizations follow, for the correlations with $\BH$ having the R.e.e.f., from Theorem 2. 
\smallskip
\par
\noindent {\bf Corollary 2}. {\it Let the correlation} $C_{f,g}(N,a)$ {\it satisfy } $\BH$. {\it Then}
$$
\supporto(C_{f,g}'(N,\cdot))\subseteq (Q), 
\enspace \hbox{\it for\enspace some } \enspace Q\in \N
$$
\par
\noindent
{\it and } 
$$
\lim_P \sum_{{r\in )P(}\atop {r>1}}{{C_{f,g}'(N,dr)}\over r}=\0(d)
$$
\par
\noindent
{\it are properties both equivalent to the {\stampatello R.e.e.f.}, of $C_{f,g}(N,a)$.} 
\smallskip
\par
\noindent {\bf Proof}. Straightforward, from Theorem 2 for $F(a)=C_{f,g}(N,a)$.\qed

\bigskip
\bigskip


\par
\noindent{\bf 1.3. Another conditional Proof of Hardy-Littlewood Conjecture, under Wintner Assumption}
\bigskip
\par
\noindent
The classical \lq \lq {\stampatello Hardy-Littlewood Conjecture}\rq \rq, for $2k-$twin primes, is the asymptotic given, once fixed an even number $2k$ ($k\ge 1$), for the {\stampatello autocorrelation} of von Mangoldt function $\Lambda$ (see [T]) of {\stampatello shift} $2k$, namely \enspace $C_{\Lambda,\Lambda}(N,2k)$ \enspace (compare Conjecture B, page 42, in [HL]) :  
$$
C_{\Lambda,\Lambda}(N,2k)\sim \SingSer(2k)N,
\quad 
\hbox{\it as}\enspace N\to \infty, 
\leqno{\HL}
$$
\par
\noindent
where the {\stampatello classical Singular Series} is defined as: 
$$
\SingSer(2k)\defineq \sum_{q=1}^{\infty}{{\mu^2(q)}\over {\varphi^2(q)}}c_q(2k)
 =2\prod_{p|2k,p>2}\left(1+{1\over {p-1}}\right)\prod_{p\; \nondivide \; 2k}\left(1-{1\over {(p-1)^2}}\right),
\quad
\forall k\in \N. 
$$ 
\par
\noindent
As a very simple consequence of our Corollary 1, we get the following result, whose Proof we sketch here, closely following the Proof of Corollary 2 in [C0]. 
\par
However, as we did in [C0], we first have to, say, \lq \lq truncate\rq \rq, the function \enspace $g(m)=\sum_{d|m}g'(d)$ \enspace (by $g'$ definition),  with the $N-$truncated divisor sum called \enspace $g_N(m)\defineq \sum_{d|m,d\le N}g'(d)$ \enspace because, then, for the correlation $C_{f,g_N}(N,a)$,  \enspace $\BH$\enspace holds (but not for \enspace $C_{f,g}(N,a)$, \enspace in general); we apply this to $f=g=\Lambda$, getting $g_N(m)=\Lambda_N(m)\defineq -\sum_{d|m,d\le N}\mu(d)\log d$, but for general $f,g:\N \rightarrow \C$ the equation $(1)$ in [C0] entails 
$$
C_{f,g}(N,a)=C_{f,g_N}(N,a)+O\Big(a\cdot \max_{n\le N}|f(n)|\cdot \max_{N<q\le N+a}|g'(q)|\Big), 
\quad
\forall a\in \N, 
$$
\par
\noindent
whence in particular
$$
C_{\Lambda,\Lambda}(N,a)=C_{\Lambda,\Lambda_N}(N,a)+O\left(a\thinspace (\log N) \thinspace (\log (N+a))\right), 
\quad
\forall a\in \N. 
\leqno{({\rm T})}
$$
\par				
\noindent
We'll use hereafter the $O-$notation of Landau [D], equivalent to Vinogradov's (in fact, $A=O(B)$ amounts to \thinspace $A\ll B$, same for \enspace $A=O_{\varepsilon}(B)$ \enspace and \enspace $A\ll_{\varepsilon} B$). 
\smallskip
\par
In fact, actually, Wintner Assumption $\WA$, instead of $\DH$, suffices to prove even more than $\HL$. 
\smallskip
\par
\noindent {\bf Corollary 3}. {\it Assuming } $\WA$ {\it for } $C_{\Lambda,\Lambda_N}(N,a)$, {\it i.e.,} 
$$
\sum_{d=1}^{\infty}{1\over d}\left| C_{\Lambda,\Lambda_N}'(N,d)\right|\, <\infty, 
$$
\par
\noindent
{\it we get a kind of Hardy-Littlewood asymptotic formula, with an absolute constant } $c>0$, {\it once } $k\in \N$ {\it is fixed} 
$$
C_{\Lambda,\Lambda}(N,2k)=\SingSer(2k)N+O\left(N\; e^{-c\sqrt{\log N}}\right).
$$
\smallskip
\par
\noindent {\bf Proof.(Sketch)} We first get the $\REEF$ for $C_{\Lambda,\Lambda_N}(N,a)$, from Corollary 1 above. Then, $({\rm T})$ above, say \lq \lq Truncation Formula\rq \rq, reduces a known calculation, performed in [C0] Corollary 2 Proof, for the RHS (Right Hand Side) of the $\REEF$ to the RHS above.\qed 

\bigskip
\bigskip

\par
\noindent{\bf 2. Lemmata for the Theorems}
\bigskip
\par
\noindent
We recall hereafter $\1_{\wp}\defineq 1$ if and only if the property $\wp$ is true and otherwise $\defineq 0$, in the formula (compare [D] and [M]): $\forall n\in \N$, $\forall a\in \Z$ 
$$
\sum_{q|n}c_q(a)=\sum_{q|n}\sum_{{d|q}\atop {d|a}}d\mu\left({q\over d}\right)
 =\sum_{{d|a}\atop {d|n}}d\sum_{{q|n}\atop {q\equiv 0\bmod d}}\mu\left({q\over d}\right)
 =\sum_{{d|a}\atop {d|n}}d\sum_{K\left|{n\over d}\right.}\mu(K)
 = \1_{n|a}\cdot n, 
\leqno{(2)}
$$
\par
\noindent
from {\stampatello Kluyver's formula}: ${\displaystyle \sum_{d|a,d|q}d\left({q\over d}\right)=c_q(a) }$ [K] and {\stampatello M\"{o}bius inversion}: \thinspace ${\displaystyle \sum_{K|m} \mu(K)=\1_{m=1} }$ \thinspace \thinspace (see [T]). 
\medskip
\par
Next elementary Lemma is most of our first Proof of Theorem 1, see $\S1.1$ above. Recall $p-$adic valuation: as usual, $v_p(a)\defineq \max\{ K\in \N_0 : p^K|a\}$, where $\N_0\defineq \N \cup \{0\}$. Recall also: $\pi(x)\defineq \left|\left\{p\in \Primes\,:\, p\le x\right\}\right|$. 
\smallskip
\par
\noindent {\bf Lemma 1}. ({\stampatello Elementary properties of smooth divisors})
\par
\noindent
{\it Let } $P$ {\it be a prime number and $F:\N \rightarrow \C$ be any arithmetic function, with Eratosthenes transform } $F'$ ({\it recall, } $F'\defineq F\ast \mu$). {\it Then, } $\forall a\in \N$ {\stampatello FIXED}, 
\medskip
\par
\noindent
$$
P\ge a 
\enspace \Rightarrow \enspace 
F(a)=\sum_{d\in (P)}{{F'(d)}\over d}\sum_{{q\in (P)}\atop {q|d}}c_q(a); 
\leqno{(3)}
$$
\par
\noindent
$$
\forall d\in \N,\enspace 
\sum_{{q\in (P)}\atop {q|d}}|c_q(a)|\le \sum_{q\in (P)}|c_q(a)|
 \le 2^{\pi(P)}a<\infty; 
\leqno{(4)}
$$
\par
\noindent
$$
\forall d\in \N,\thinspace 
\sum_{{q\in (P)}\atop {q|d}}c_q(a)=\sum_{q\left|d_{(P)}\right.}c_q(a)
 =\1_{d_{(P)}\left|a\right.}\cdot d_{(P)},\thinspace \hbox{\it with } d_{(P)}\defineq \prod_{p\le P} p^{v_p(d)},\thinspace \hbox{\it whence}\enspace 0\le 
\sum_{{q\in (P)}\atop {q|d}}c_q(a)\le a. 
\leqno{(5)}
$$
\medskip
\par
\noindent {\bf Proof}. Fix $a\in \N$, take $P\ge a$, $P\in \Primes$, write $\1_{d|a}$ from $(2)$, getting 
$$
F(a)=\sum_{{d\in (P)}\atop {d|a}}F'(d)
 =\sum_{d\in (P)}{{F'(d)}\over d}\sum_{q|d}c_q(a), 
  =\sum_{d\in (P)}{{F'(d)}\over d}\sum_{{q\in (P)}\atop {q|d}}c_q(a), 
$$
\par				
\noindent
because $(P)$ is {\it divisor-closed}, namely $d\in (P)$ and $q|d$ imply $q\in (P)$, giving $(3)$.\hfill QED
\par
\noindent
Since $|c_q(a)|$ is a multiplicative function of $q\in\N$, see Fact 1 \& Main Lemma of [C2] for details, 
$$
\sum_{q\in (P)}|c_q(a)|=\prod_{p\le P}\left(\sum_{K=0}^{v_p(a)}\varphi(p^K)+p^{v_p(a)}\right)
 =\prod_{p\le P}\left(2p^{v_p(a)}\right)\le 2^{\pi(P)} a, 
$$
\par
\noindent
providing $(4)$.\hfill QED 
\par
\noindent
The condition \lq \lq $q\in (P)$ and $q|d$\rq \rq, by definition of $d_{(P)}$, is equivalent to the single condition $q$ divides $d_{(P)}$, so $(2)$ with $n=d_{(P)}$ entails $(5)$.\hfill QED
\par
\noindent
The Lemma is completely settled.\qed

\medskip

\par
See that the main reason why our Theorem 1 works for Ramanujan smooth expansions but not for Ramanujan expansions is inside property $(5)$ above; in fact, if we wish, say, to get the same $d-$independent bound for usual partial sums, we should consider (as Wintner does explicitly, see [W] page 31) 
$$
\sum_{{q\le Q}\atop {q|d}}c_q(a), 
$$
\par
\noindent
which has not a closed expression similar to the one in $(5)$ : this time, the multiplicative structure is, say, broken by the interval constraint. 
\par
\noindent
Note the, say, very simple structure of Lemma 1: once added $\WA$, the proof of Theorem 1 is immediate.

\medskip

\par
We wish to prove a kind of equivalence condition, for the convergence for Ramanujan smooth series with Wintner coefficients (see Lemma 3, next section). So, in next Lemma we, say, decompose in a regular part (over $P-$smooth numbers), containing Wintner coefficients, and an irregular part (over $P-$sifted numbers), containing Eratosthenes transform. We do apply this decomposition in $\S3$ : an alternative (much) shorter Proof of Theorem 1, then, is immediate. We might say that Wintner's Dream Theorem is a straightforward application of Wintner's ($P-$)Orthogonal Decomposition, i.e., next Lemma 2. 

\bigskip

\par
\noindent
From M\"obius inversion [T] quoted above, abbreviating $P\P \defineq {\displaystyle \prod_{p\le P}p }$, we get the useful formul\ae: 
$$
\1_{(a,b)=1} = \sum_{{K|a}\atop {K|b}}\mu(K),
\enspace \forall a,b\in \N 
\quad \Rightarrow \quad 
\1_{r\in )P(} = \sum_{{K|r}\atop {K\left|P\P\right.}}\mu(K)
 =\sum_{{K|r}\atop {K\in (P)}}\mu(K). 
\leqno{(6)}
$$
\par
A kind of \lq \lq arithmetic orthogonality among indices\rq \rq, say, allows to decompose $F'$ in $(7)$, then $F$ in $(8)$: both vital, for our arguments. 
\smallskip
\par
\noindent {\bf Lemma 2}. ({\stampatello Wintner orthogonal decomposition})
\par
\noindent
{\it Let } $F:\N\rightarrow\C$ {\it have all Wintner coefficients, say } $\exists \WintnerT F \, :\, \N \rightarrow \C .$ {\it Then} 
$$
\forall d\in \N, \enspace \forall P\in \Primes, 
\quad 
F'(d)=d\sum_{K\in (P)}\mu(K)\left(\Wintner_{dK} F\right)-\sum_{{r\in )P(}\atop {r>1}}{{F'(dr)}\over r},
\leqno{(7)}
$$
\par
\noindent
{\it whence}
$$
\forall d\in \N, \quad 
F'(d)=\lim_P \left(d\sum_{K\in (P)}\mu(K)\left(\Wintner_{dK} F\right)-\sum_{{r\in )P(}\atop {r>1}}{{F'(dr)}\over r}\right). 
$$
\par
\noindent
{\it If we join the hypothesis: } $\WintnerT F$ {\it smooth-supported, say } $\supporto(\WintnerT F)\subseteq (Q)$, {\it we get} 
$$
\forall d\in \N, \quad 
F'(d)=\1_{d\in (Q)}\cdot d\cdot \sum_{K\in (Q)}\mu(K)\left(\Wintner_{dK} F\right)-\lim_P \sum_{{r\in )P(}\atop {r>1}}{{F'(dr)}\over r}, 
$$
\par				
\noindent
{\it whence in particular for finite support, say, } $\supporto(\WintnerT F)\subseteq [1,Q]$, {\it this entails} 
$$
\forall d\in \N, \quad 
F'(d)=d\sum_{K\le {Q\over d}}\mu(K)\left(\Wintner_{dK} F\right)-\lim_P \sum_{{r\in )P(}\atop {r>1}}{{F'(dr)}\over r}. 
$$
\par
\noindent
{\it Summing $(7)$ over the divisors $d$ of $a$, we obtain } ({\it however } $P\in \Primes$, {\it here}) 
$$
\forall a\in \N, \enspace \forall P\ge a, 
\quad 
F(a)=\sum_{q\in (P)}\left(\Wintner_q F\right)c_q(a)-\sum_{d|a}\sum_{{r\in )P(}\atop {r>1}}{{F'(dr)}\over r},
\leqno{(8)}
$$
\par
\noindent
{\it whence}
$$
\forall a\in \N, \enspace \quad 
F(a)=\lim_P\left(\sum_{q\in (P)}\left(\Wintner_q F\right)c_q(a)-\sum_{d|a}\sum_{{r\in )P(}\atop {r>1}}{{F'(dr)}\over r}\right). 
$$
\par
\noindent
{\it This time, } $\supporto(\WintnerT F)\subseteq (Q)$ {\it gives}
$$
\forall a\in \N, \enspace \quad 
F(a)=\sum_{q\in (Q)}\left(\Wintner_q F\right)c_q(a)-\lim_P \sum_{d|a}\sum_{{r\in )P(}\atop {r>1}}{{F'(dr)}\over r}, 
$$
\par
\noindent
{\it in particular } $\supporto(\WintnerT F)\subseteq [1,Q]$ {\it entails} 
$$
\forall a\in \N, \enspace \quad 
F(a)=\sum_{q\le Q}\left(\Wintner_q F\right)c_q(a)-\lim_P \sum_{d|a}\sum_{{r\in )P(}\atop {r>1}}{{F'(dr)}\over r}. 
$$
\medskip
\par
\noindent {\bf Remark 5.} The $r-$series above (defined in Remark 3, $\S1$) is called the {\stampatello Irregular series}, $\IrrPdF$, {\it of $\underline{\hbox{\it argument}}$ $d\in \N$, $\underline{\hbox{\it over}}$ the prime $P\in \Primes$, $\underline{\hbox{\it relative}}$ to $F:\N \rightarrow \C$}, and the following Proof implies it converges in $\C$, when \thinspace $\WintnerT F$ exists.\hfill $\diamond$ 
\medskip
\par
\noindent {\bf Proof}. In order to prove $(7)$, we fix $d\in \N$ and $P\in \Primes$, considering  
$$
\sum_{r\in )P(}{{F'(dr)}\over r}=\lim_x \sum_{{r\in )P(}\atop {r\le x}}{{F'(dr)}\over r}
 =\lim_x \sum_{r\le x}{{F'(dr)}\over r}\sum_{{K|r}\atop {K\in (P)}}\mu(K)
  =d\lim_x \sum_{K\in (P)}\mu(K)\sum_{{r\le x}\atop {r\equiv 0\bmod K}}{{F'(dr)}\over {dr}}, 
$$
\par
\noindent
thanks to $(6)$; the $K-sum$, thanks to $\mu(K)$, is over the square-free $K$ and, furthermore, the condition that $K$ divides the $P-${\it primorial} (abbreviated $P\P$) amounts to $K\in (P)$, from: $K$ square-free; in all, thanks to the fact: $\mu$ is supported in square-free numbers, say \lq \lq {\it M\"{o}bius vertical limit}\rq \rq, this $K-$sum is FINITE AND clearly NOT DEPENDING ON $x$, giving : 
$$
\sum_{r\in )P(}{{F'(dr)}\over r}=d\sum_{K\in (P)}\mu(K)\lim_x \sum_{{r\le x}\atop {r\equiv 0\bmod K}}{{F'(dr)}\over {dr}}
 =d\sum_{K\in (P)}\mu(K)\left(\Wintner_{dK}\; F\right), 
$$
\par
\noindent
thanks to the definition of Wintner coefficients (all series converging for them, since $\exists \WintnerT F$). 
\par
Separating the contribute of \enspace $r=1$ \enspace in the $r-$series settles $(7)$ proof.\hfill QED
\par
\noindent
Joining \enspace $\supporto(\WintnerT F)\subseteq (Q)$, whenever $P\ge Q$, then $(7)$, in particular for the case $\supporto(\WintnerT F)\subseteq [1,Q]$, entails both the two particular formul\ae, after $(7)$.\hfill QED 
\par
Then, $(8)$ comes from $(7)$ summing over $d|a$, with Kluyver formula and: $d|a,P\ge a$ $\Rightarrow$ $d\in (P)$.\hfill QED 
\par
\noindent
The two particular formul\ae \thinspace after $(8)$ follow from $(8)$, as we saw for $(7)$, above.\qed 
\medskip
\par
\noindent {\bf Remark 6.} As it's clear from the Proof, in case \enspace $\supporto(\WintnerT F)\subseteq (Q)$ (in particular, whenever we have $\supporto(\WintnerT F)\subseteq [1,Q]$, too) we get that {\stampatello the Irregular series} (defined in Remark 3), $\IrrPdF$, {\stampatello is constant $\forall P\ge Q$, w.r.t. the prime $P$, uniformly in the argument $d\in \N$} 
$$
\sum_{{r\in )P(}\atop {r>1}}{{F'(dr)}\over r}=\sum_{{r\in )Q(}\atop {r>1}}{{F'(dr)}\over r}, 
\enspace \forall P>Q
\quad \Rightarrow \quad
\lim_P \sum_{{r\in )P(}\atop {r>1}}{{F'(dr)}\over r} = \sum_{{r\in )Q(}\atop {r>1}}{{F'(dr)}\over r}, 
\leqno{(9)}
$$
\par
\noindent
i.e., the LHS (Left Hand Side) of $(9)$, as a function of $P\in\Primes$, is constant $\forall P\ge Q$, uniformly $\forall d\in \N$. Then, \hfill
\par
\noindent
notice that (assuming, as we can, that our $Q$ is prime) when $\underline{\hbox{\rm a\enspace fortiori}}$ $\supporto(F')\subseteq (Q)$, \enspace we have from $(9)$: \enspace $\lim_P \IrrPdF = {\rm Irr}^{(Q)}_d\,F = \0(d)$, because $F'(dr) = \0(d)$, $\forall r\in )Q(\,\backslash \{1\}$.\hfill $\diamond$ 

\vfill
\eject

\par				
\noindent{\bf 3. A deeper look into Ramanujan smooth expansions: Ramanujan-Wintner smooth expansions}
\bigskip
\par
\noindent
A more careful analysis yields in fact the more general result, for THE Ramanujan smooth expansion, with Wintner coefficients. We write THE to highlight its uniqueness, clear from the choice of coefficients $G(q):=\Wintner_q F$. As we'll see in section 5.1, once fixed $F$ (esp., $F=\0$), we may have many $G$, in a Ramanujan smooth expansion. 
\smallskip
\par
\noindent {\bf Lemma 3}. ({\stampatello characterizing $F$ having Ramanujan-Wintner smooth expansion})
\par
\noindent
{\it Let } $F:\N\rightarrow\C$ {\it have all the Wintner coefficients. Then, $\forall d\in \N$ fixed}
$$
F'(d)=\lim_{P}\; d\sum_{K\in (P)}\mu(K)\left( \Wintner_{dK} F \right)
\enspace \Longleftrightarrow \enspace 
\lim_{P} \sum_{{r\in )P(}\atop {r>1}}{{F'(dr)}\over r} = 0, 
$$
{\it whence, $\forall a\in \N$ fixed, } 
$$
F(a)=\lim_{P} \sum_{q\in (P)}\left( \Wintner_q F \right)c_q(a)
\enspace \Longleftrightarrow \enspace 
\lim_{P} \sum_{d|a}\sum_{{r\in )P(}\atop {r>1}}{{F'(dr)}\over r} = 0. 
$$
\par
\noindent {\bf Proof}. From Lemma 2, passing to the limit over $P\in \Primes$, we get first equivalence from $(7)$ and second one from $(8)$.\qed

\bigskip

\par
\noindent
We give an easy property (next Proposition), connecting $|F'|\ast \1$ to $F$. 

\medskip

From above Lemma 3 and the trivial implication, $\forall d\in \N$, 
$$
\lim_P \sum_{{r\in )P(}\atop {r>1}}{{|F'(dr)|}\over r} = 0
\quad \Rightarrow \quad 
\lim_P \sum_{{r\in )P(}\atop {r>1}}{{F'(dr)}\over r} = 0
$$
\par
\noindent
we easily prove the following. We abbreviate \lq \lq RWE\rq \rq, for \lq \lq Ramanujan-Wintner expansion\rq \rq: Ramanujan expansion with Wintner coefficients. Joining \lq \lq {\it smooth}\rq \rq, hereafter, amounts, as above, to requiring {\it smooth partial sums}. 
\smallskip
\par
\noindent {\bf Proposition 1}. {\it Given any $F:\N \rightarrow \C$, we have} 
$$
|F'|\ast \1 \enspace \hbox{\it has\enspace smooth\enspace {\rm RWE}} 
\quad \Rightarrow \quad 
F \enspace \hbox{\it has\enspace smooth\enspace {\rm RWE}}. 
$$ 

\medskip

\par
\noindent
Notice that, actually, this can also be proved following Theorem 1 proof in $\S1$. 

\medskip

\par
\noindent
By the way, we give now a shorter proof of this Theorem. 

\medskip

\par
\noindent
\centerline{\bf Alternative proof of Wintner's Dream Theorem}
\bigskip 
\par
\noindent {\bf Proof}. Using $(8)$ of Lemma 2, i.e., applying Wintner Orthogonal Decomposition to $F$, it suffices to prove:
$$
\WA
\enspace \Rightarrow \enspace 
\lim_P \sum_{{r\in )P(}\atop{r>1}}{{F'(dr)}\over r}=\0(d), 
$$
\par
\noindent
so, fix $d\in \N$ and consider:\enspace ${\displaystyle \sum_{{r\in )P(}\atop{r>1}}{{|F'(dr)|}\over {dr}}\le \sum_{{m>dP}\atop{m\equiv 0\bmod d}}{{|F'(m)|}\over m}\le \sum_{n>P}{{|F'(n)|}\over n} }$ is infinitesimal, as\enspace $P\to \infty$.\qed 

\vfill
\eject

\par				
\noindent{\bf 4. Smooth coefficients in Ramanujan expansions}
\bigskip
\par
\noindent
In our \lq \lq {\it A smooth shift approach for a Ramanujan expansion}\rq \rq, [C1], we introduced the {\stampatello smooth restriction}, to $P-$smooth numbers (here $P\in \Primes$ is fixed), {\stampatello of any given arithmetic function} $F:\N \rightarrow \C$, 
$$
F_{(P)}(a)\defineq \sum_{{d\in (P)}\atop {d|a}}F'(d), 
\quad 
\forall a\in \N, 
$$
\par
\noindent
which is, so to speak, the origin of {\stampatello Carmichael's \& Wintner's} \lq \lq $P-${\stampatello smooth coefficients}\rq \rq: 
$$
\Carmichael_q^{(P)} F \defineq {1\over {\varphi(q)}}\lim_{x\to \infty}{1\over x}\sum_{a\le x}\left(\sum_{{d\in (P)}\atop {d|a}}F'(d)\right)c_q(a)
 = \Carmichael_q F_{(P)}, 
\qquad
\forall q\in \N 
$$
\par
\noindent 
and 
$$
\Wintner_q^{(P)} F \defineq \sum_{{d\in (P)}\atop {d\equiv 0\bmod q}}{{F'(d)}\over d}
 = \Wintner_q F_{(P)}, 
\qquad
\forall q\in \N.
$$
\par
\noindent 
We say in the following that they exist, whenever the relative limits exist in $\C$ (for Wintner's, the limit of partial sums, i.e., the series converges in $\C$). The existence of all coefficients (for $P$ fixed), $\forall q\in \N$, is expressed saying : $\exists \CarmichaelTP F$ or, resp., $\exists \WintnerTP F$. They are, resp., the {\it Carmichael $P-$smooth transform} and {\it Wintner $P-$smooth transform}; when they both exist and are the same, we indicate them as $\widehat{F}^{(P)}$, say the {\it Carmichael-Wintner $P-$smooth transform}. In [C1] we gave Theorem 1 \& Corollary 1, which we will generalize here to the following Theorem 1'[C1] \& Corollary 1'[C1]. 
\par 
The interest in these $P-$smooth transforms comes from the following Theorem 1'[C1], giving a kind of \lq \lq limit expansion\rq \rq, for any \lq \lq reasonable\rq \rq, say, arithmetic function $F$ : we mean that $\exists \WintnerTP F$, $\forall P\in \Primes$. 
\medskip
\par 
It's based on next elementary Lemma, an immediate application of \lq \lq {\stampatello Ramanujan vertical limit}\rq \rq:
$$
c_q(a)\neq 0
\enspace \Rightarrow \enspace
v_p(q)\le v_p(a)+1,\enspace \forall p|q.
\leqno{\Rvl}
$$ 
\smallskip
\par
\noindent {\bf Lemma 4}. {\it Let } $F:\N\rightarrow\C$ {\it have all the $P-$smooth $q-$th Wintner coefficients}:  $\exists \WintnerPqF$, $\forall P\in \Primes, \forall q \in \N$. {\it Then}
$$
\forall a\in \N, \forall P\in \Primes, P\ge a,  
\thinspace \enspace \thinspace  
F(a)=\sum_{q\in (P)}\left(\WintnerPqF \right)c_q(a), 
$$
\par\noindent{\it whence} \quad
$$
\forall a\in \N, 
\quad 
F(a)=\lim_P \sum_{q\in (P)}\left(\WintnerPqF \right)c_q(a). 
$$ 
\smallskip
\par
\noindent {\bf Proof}. Fix $a\in \N$ and choose a prime $P\ge a$ so that $d|a$ $\Rightarrow $ $d\in (P)$ and from $(2)$ 
$$
F(a)=\lim_x \sum_{{d|a}\atop {d\in (P)\, ,\, d\le x}}F'(d)
 =\lim_x \sum_{{d\in (P)}\atop {d\le x}}{{F'(d)}\over d}\sum_{q|d}c_q(a), 
$$
\par
\noindent
whence $d\in (P),q|d$ $\Rightarrow $ $q\in (P)$ proves, from $\Rvl$, that the sum over $q$ is both finite (in terms of $a\in \N$ and $P\ge a$, $P\in \Primes$) and is not dependent on $x$ (going to infinity); giving 
$$
F(a)=\sum_{q\in (P)}c_q(a)\lim_x \sum_{{d\in (P),d\le x}\atop {d\equiv 0\bmod q}}{{F'(d)}\over d}
 =\sum_{q\in (P)}\left(\WintnerPqF \right)c_q(a), 
$$
\par
\noindent
from the definition of $P-$smooth $q-$th Wintner coefficient.\qed 
\smallskip
\par
\noindent {\bf Remark 7.} See that, usually, the exchange of two summations, typically over $d$ and $q$ like in the above, needs a double series (over $d,q$) {\sl absolute convergence}, while here the $\Rvl$ property allows weaker hypotheses. Also, notice that (whatever $P\in \Primes$ is fixed) the condition: $\exists \WintnerPqF$, $\forall q\not\in (P)$, {\bf isn't strictly required}.\hfill $\diamond$ 

\vfill
\eject

\par				
This Lemma is {\sl very powerful} : each time we have hypotheses ensuring the existence of $\underline{\hbox{\rm all}}$ the $P-${\stampatello smooth $q-$th Wintner coefficients}, we get a kind of \lq \lq {\stampatello Ramanujan-Wintner local expansion}\rq \rq (with a {\stampatello smooth summation} of partial sums \& $P-${\stampatello smooth Wintner coefficients}). The only problem is the, say, \lq \lq {\stampatello local nature of coefficients}\rq \rq, that usually are unkown; while, of course, {\stampatello Wintner coefficients} have better chances to be easily calculated: for example, under suitable hypotheses, they are exactly the {\stampatello Carmichael coefficients}. This happens under $\WA$ above, as proved by Wintner (see [C3]). 

\bigskip

\par
The same $\WA$ ensures that 
$$
\lim_P \WintnerPqF = \Wintner_q F,
\qquad
\hbox{\stampatello uniformly}\enspace \forall q\in \N, 
$$
\par
\noindent
thanks to the absolute convergence for the series inside $\WA$. 
\par
Wintner Assumption, actually, suffices (see Theorem 1) for $F:\N \rightarrow \C$ to get $\underline{\hbox{\stampatello the}}$ {\stampatello Ramanujan-Wintner Smooth expansion} for $F$. 
\par
Can we get the same expansion under a weaker hypothesis ? Well, our Theorem 1 proof reveals this \lq \lq at once\rq \rq. \enspace From the point of view of  Wintner coefficients, next result is, in fact, a generalization of our Theorem 1 above. \hfil We give it here, as its hypotheses are a bit more technical than Theorem 1 ones. 
\medskip 

\par
{\stampatello Wintner's Smooth Assumption} $\WSA$, following, is a less general constraint than $\WA$ above: 
$$
\lim_P \sum_{d\not\in (P)}{{|F'(d)|}\over d}=0
\leqno{\WSA}
$$
\par
\noindent
where we implicitly agree that: given our $F:\N \rightarrow \C$, there exists a prime $P_F$ (depending ONLY on $F$), such that each series above over $d\not\in (P)$ converges $\forall P>P_F$ and, then, above limit over $P$ exists and vanishes. This $\WSA$ alone proves that $F$ converges, with smooth partial summations, to {\sl its} Ramanujan-Wintner Smooth expansion and this is already proved, in Theorem 1 proof (see its end) !
\par
Joining two technical hypotheses about (classic \& smooth) Wintner coefficients we also get, say for free, other two informations: see next result. 
\smallskip
\par
\noindent {\bf Theorem 1(Smooth version)}. {\it Let } $F:\N\rightarrow\C$ {\it have all the $q-$th Wintner coefficients (i.e., $\exists \WintnerT F$) and all the $P-$smooth $q-$th Wintner coefficients (i.e., $\exists \WintnerTP F, \forall P$), } $\forall P\in \Primes, \, \forall q \in \N$. {\it Assume } $\WSA$.  {\it Then}
$$
\lim_P \WintnerPqF = \Wintner_q F,
\qquad
\hbox{\stampatello uniformly}\enspace \forall q\in \N 
\leqno{(\ast)}
$$ 
\par
\noindent 
{\it and} 
$$
\lim_P \IrrPdF = \0(d), 
\qquad
\hbox{\stampatello pointwisely}\enspace \forall d\in \N, 
\leqno{(\ast\ast)}
$$ 
{\it whence}
$$
F(a)=\lim_P \sum_{q\in (P)}\left(\Wintner_q F\right)c_q(a), 
\qquad
\hbox{\stampatello pointwisely}\enspace \forall a\in \N. 
\leqno{(\ast\ast\ast)}
$$ 
\smallskip
\par
\noindent {\bf Proof}. Above $(\ast)$ follows immediately from
$$
\left| \Wintner_q F - \WintnerPqF \right|\le \sum_{{d\not\in (P)}\atop {d\equiv 0\bmod q}}{{|F'(d)|}\over d}
 \le \sum_{d\not\in (P)}{{|F'(d)|}\over d}; 
$$
\par
\noindent 
for $(\ast\ast)$ use Lemma 2 to prove the convergence, whence existence of $\IrrPdF$ for all $P>P_F$ and $\forall d\in \N$, then $\forall d\in \N$ FIXED
$$
\left| \IrrPdF \right|\le \sum_{{r\in )P(}\atop {r>1}}{{\left|F'(dr)\right|}\over r} 
\le d\sum_{m\not \in (P)}{{\left|F'(m)\right|}\over m}, 
$$
\par
\noindent 
as $m=dr$ with $r\in )P($ and $r>1$ imply $\exists p>P$ : $p|m$ $\Rightarrow$ $m\not \in (P)$; for $(\ast\ast\ast)$ use $(\ast\ast)$ just proved and the characterization of Lemma 3.\qed 

\vfill
\eject

\par				
\noindent
See that previous result, in particular, needs the existence of \enspace $\WintnerTP F$, $\forall P\in \Primes$. 

\par
Theorem 1 [C1] (version 3) ensures this existence when $F$ satisfies the {\stampatello Ramanujan Conjecture} and has {\stampatello only smooth divisors}. Our Lemma 4 above allows to generalize this to next result, Theorem 1'[C1], in which the {\stampatello Ramanujan Conjecture} alone implies: $\exists \WintnerTP F$, $\forall P\in \Primes$. Next Theorem 1'[C1] has [C1] to distinguish from present Theorem 1. 
\smallskip
\par
\noindent {\bf Theorem 1'[C1]}. {\it Let } $F:\N\rightarrow\C$ {\it satisfy Ramanujan Conjecture. Then, say,} 
$$
\widehat{F}^{(P)}(q) \defineq \WintnerPqF 
 = \CarmichaelPqF 
  = \prod_{p\le P}\left(1-{1\over p}\right){1\over {\varphi(q)}}\sum_{t\in (P)}{{F(t)}\over t}c_q(t),
\quad \forall P\in \Primes, \, \forall q\in (P)
$$ 
{\it and}
$$
\forall a\in \N, \forall P\in \Primes, P\ge a, 
F(a)=\sum_{q\in (P)}\widehat{F}^{(P)}(q)c_q(a),
$$ 
\par
\noindent
{\it whence} 
$$
\forall a\in \N, 
F(a)=\lim_P \sum_{q\in (P)}\widehat{F}^{(P)}(q)c_q(a). 
$$ 
\par
\noindent
{\it In particular, all $F:\N \rightarrow \C$ satisfying Ramanujan Conjecture are pointwise limits, over primes $P\to \infty$, of \lq \lq finite Ramanujan expansions\rq \rq, with \lq \lq Carmichael-Wintner\rq \rq \enspace $P-$smooth coefficients }({\it i.e., $\widehat{F}^{(P)}(q)$ above}). 
\smallskip
\par
\noindent {\bf Remark 8.} The main \lq \lq defect\rq \rq, so to speak, is the fact that the coefficients may change, as $P$ changes.\hfill $\diamond$ 
\smallskip
\par
\noindent {\bf Proof}. First of all, the explicit formula above for $\CarmichaelPqF$ is proved in $(ii)$ of Theorem 1 [C1]; then, the coincidence, for all $P\in \Primes$, of $\CarmichaelTP F$ and $\WintnerTP F$ was proved in Th.m 1, $(i)$ [C1] (in which these coefficients were born). The formula for $F$, then, was proved in [C1] (Theorem 1 proof), in case $F'$ is supported over $P-$smooth numbers: this is implicit here, assuming $P\ge a$. So, present second part is more general than Theorem 1 in [C1].\qed

\medskip

Its immediate application follows, to the Correlations. Again, we join [C1] to distinguish from present Corollary 1. 
\smallskip
\par
\noindent {\bf Corollary 1'[C1]}. {\it Let } $C_{f,g}(N,a)$, {\it the correlation of any couple} $f,g: \N \rightarrow \C$, {\it satisfy Ramanujan Conjecture}. {\it Then} 
$$
\forall a\in \N, 
\quad 
C_{f,g}(N,a)=\lim_P \sum_{q\in (P)}\widehat{C_{f,g}}^{(P)}(N,q)c_q(a), 
$$ 
\par
\noindent
{\it where the \enspace $\widehat{C_{f,g}}^{(P)}(N,q)$ \enspace are \lq \lq Carmichael-Wintner $P-$smooth $q-$th coefficients\rq \rq, for correlations, see } [C1]. 
\par
\noindent {\bf Proof}. Apply Theorem 1'[C1] to $F(a)=C_{f,g}(N,a)$.\qed
\medskip
\par
\noindent {\bf Remark 9.} For correlations satisfying $\BH$, Ramanujan Conjecture follows: see $(ii)$,Proposition 1 [C1].\hfill $\diamond$ 

\medskip

\par
More in general, Ramanujan Conjecture for $F$ is not required, if we wish to get the existence of all $\WintnerTP F$ ($\forall P\in \Primes$). In fact, it holds also for the $F:\N \rightarrow \C$ satisfying 
$$
\exists \delta<1 \, : \, F(a)\ll_{\delta} a^{\delta},
\quad 
\hbox{\rm as }\enspace a\to \infty. 
\leqno{\NSL}
$$
\par
\noindent
This is, by M\"{o}bius inversion [T], equivalent to the same for $F'$, the {\it Eratosthenes Transform} of our $F$. 
\smallskip
\par
This property of {\stampatello Neat Sub-Linearity}, actually, implies even more than the existence of all $\WintnerPqF$: 
$$
\NSL
\enspace \Rightarrow \enspace 
\sum_{{d\in (P)}\atop {d\equiv 0\bmod q}}{{|F'(d)|}\over d}\ll_{\delta} \sum_{d\in (P)}d^{\delta-1}
 \ll_{\delta,P} 1,
\enspace 
\hbox{\stampatello uniformly}\enspace \forall q\in \N. 
$$
\par
\noindent
This is another application of $(1)$ above. 

\medskip

\par
However, the existence of all $\WintnerPqF$ also follows from another hypothesis for $F$: 
$$
F'=\mu^2 \cdot F', 
\leqno{\IPP}
$$
\par				
\noindent
i.e. $F'$ {\sl is \enspace square-free\enspace supported}, say, $F$ \lq \lq {\stampatello Ignores Prime-Powers}\rq \rq. Equivalently, $F(a)$ depends ONLY on $\kappa(a)\defineq \prod_{p|a}p$ (with $\kappa(1)\defineq 1$ for the void product), the {\stampatello square-free kernel} of our $a\in \N$: we express this as $F=F\circ \kappa$ (with \lq \lq $\circ$\rq \rq, here, the usual composition of functions), namely $F(a)=F(\kappa(a)), \forall a\in \N$. 
\par
In fact, $\IPP$ implies that 
$$
\forall q\in \N,
\enspace
\sum_{{d\le x}\atop {{d\in (P)}\atop {d\equiv 0\bmod q}}}{{F'(d)}\over d} = \sum_{{d\le x}\atop {{d\in (P)}\atop {d\equiv 0\bmod q}}}{{\mu^2(d)F'(d)}\over d}
$$
\par
\noindent
has a finite limit in complex numbers, as $x\to \infty$, since previous summation's support is bounded, having cardinality bounded uniformly $\forall q\in \N$ as 
$$
\left|\left\{d\in (P) : \mu^2(d)=1\right\}\right|=2^{\pi(P)}.
$$
\par
\noindent
(Recall: $\pi(P)=$ number of primes $\le P$ and all square-free $d$ with prime factors $\le P$ are $2^{\pi(P)}$, of course.)

\medskip

\par 
For classic Carmichael \& Wintner coefficients, Wintner discovered their coincidence, whenever His $\WA$ holds. 

\par
Actually, a little bit more generally, under the following hypothesis: 
$$
\lim_x {1\over x}\sum_{d\le x}\left| F'(d)\right|=0, 
\leqno{\ETD}
$$
\par
\noindent
say, \lq \lq {\stampatello Eratosthenes Transform Decay}\rq \rq, equivalent to the vanishing of $|F'(d)|$ mean-value (esp., see [C3], Remark 7), we get\enspace $\CarmichaelT F = \WintnerT F$ again, from following Lemma (compare the proof of $(5)$ in [C3]). 
\par
See that $\WA \Rightarrow \ETD$ (from quoted proof), but the converse implication doesn't hold (esp., we may take $F'(d)=1/\log d$, $\forall d>1$). 

\medskip

\par
However, just like $\ETD$ implies coincidence of Carmichael \& Wintner coefficients, say a classic consequence, it also implies, for all fixed primes $P$, the coincidence, say, of Carmichael \& Wintner $P-$smooth transforms: $\WintnerTP F=\CarmichaelTP F$. These two consequences for $F$, under $\ETD$, hold thanks to the next Lemma. (In which we express the {\sl proximity of partial sums} up to $x\in \N$, say; in fact, the coefficients exist, by definition, if and only if the $x-$limit exists in complex numbers.) 
\smallskip
\par
\noindent {\bf Lemma 5}. ({\stampatello Links between classic \& smooth Carmichael/Wintner coefficients}). 
\par
\noindent
{\it Given any $F:\N \rightarrow \C$, } $\forall P\in \Primes$, $\forall q\in \N$, $\exists C(q)>0$ {\it such that, } $\forall x\in \N$, 
$$
\left|{1\over x}\sum_{a\le x}F(a){{c_q(a)}\over {\varphi(q)}}-\sum_{{d\le x}\atop {d\equiv 0\bmod q}}{{F'(d)}\over d}\right|\le {{C(q)}\over x}\sum_{d\le x}|F'(d)| 
$$ 
\par
\noindent
{\it and} 
$$
\left|{1\over x}\sum_{a\le x}\sum_{{d\in (P)}\atop {d|a}}F'(d){{c_q(a)}\over {\varphi(q)}}-\sum_{{d\le x}\atop {{d\in (P)}\atop {d\equiv 0\bmod q}}}{{F'(d)}\over d}\right|\le {{C(q)}\over x}\sum_{{d\le x}\atop {d\in (P)}}|F'(d)|. 
$$ 
\smallskip
\par
\noindent {\bf Remark 10.} Notice that the positive constant $C(q)$ depends $\underline{\hbox{\it only}}$ on $q\in \N$.\hfill $\diamond$ 

\medskip

\par
We briefly prove this Lemma from the following elementary \lq \lq fact\rq \rq (a kind of short Lemma). 
\smallskip
\par
\noindent {\bf Fact 1}. {\it Once fixed $d,q\in\N$, we get } ${\displaystyle \sum_{m\le {x\over d}}c_q(dm)=\1_{q|d} \cdot \varphi(q) \cdot {x\over d}+O_q(1) }$, $\forall x\in \N$. 
\par
\noindent
In fact, use ${\displaystyle c_q(dm)=\sum_{j\in \Z_q^*}e_q(jdm) }$, whence $q\not |d$ $\Rightarrow$ ${\displaystyle  \sum_{m\le {x\over d}}c_q(dm)=O\left(\sum_{{\ell|q}\atop {{\ell<q}\atop {(q,d)=\ell}}}\sum_{j\in \Z_q^*}{1\over {\left\Vert {{j(d/\ell)}\over {q/\ell}}\right\Vert}}\right) }=O_q(1)$. 
\par				
\noindent {\bf Proof(Lemma 5)}. We prove second inequality (first is similar), exchanging sums \& applying Fact 1 : 
$$
{1\over x}\sum_{a\le x}\sum_{{d\in (P)}\atop {d|a}}F'(d){{c_q(a)}\over {\varphi(q)}}={1\over x}\sum_{{d\in (P)}\atop {d\le x}}F'(d)\cdot {1\over {\varphi(q)}}\sum_{m\le {x\over d}}c_q(dm)
 = \sum_{{d\le x}\atop {{d\in (P)}\atop {d\equiv 0\bmod q}}}{{F'(d)}\over d}+O_q\left({1\over x}\sum_{{d\le x}\atop {d\in (P)}}|F'(d)|\right).
$$
\qed

\medskip

\par
From previous Lemma, we get, for $P-$smooth coefficients, using $\exists \WintnerTP F$, $\forall P\in \Primes$, as we saw, for the $F:\N \rightarrow \C$ satisfying $\NSL$ $\underline{\hbox{\rm or}}$ $\IPP$, the equation $\CarmichaelTP F=\WintnerTP F$, $\forall P\in \Primes$, too. In fact, the remainder in previous lemma goes to $0$ as $x\to \infty$ : under $\NSL$ by $(1)$, while under $\IPP$ because the $d-sum$ is bounded (w.r.t. $x$). Thus, previous Lemma 5 implies next Lemma 6. 
\smallskip
\par
\noindent {\bf Lemma 6}. ({\stampatello Two conditions for coincidence of $P-$smooth Carmichael/Wintner transforms}). 
\par
\noindent
{\it Let $F:\N \rightarrow \C$ satisfy } $\NSL$ $\underline{\hbox{\it or}}$ $\IPP$. {\it Then, } $\forall P\in \Primes$, $\CarmichaelTP F = \WintnerTP F$. 

\medskip

\par
The hypothesis $\ETD$ is $\underline{\hbox{\it only}}$ able to prove \thinspace $\CarmichaelTP F = \WintnerTP F$ \thinspace (from Lemma 5 above), when we already know that \enspace $\exists \CarmichaelTP F$ \enspace  $\underline{\hbox{\it or}}$ \enspace $\exists \WintnerTP F$. Notwithstanding its greater generality w.r.t. $\WA$, our $\ETD$ can  $\underline{\hbox{\it only}}$ prove \thinspace $\CarmichaelT F = \WintnerT F$ \thinspace (from quoted Lemma) if we know, again, that  $\underline{\hbox{\it at\enspace least\enspace one}}$ of these two transforms exists. (Of course, instead, $\WA$ $\Rightarrow$ $\exists \WintnerT F$, immediately.) 
\medskip

\par
Applying Lemma 6, we easily prove (leaving as exercises) the following two properties of uniqueness, for the $P-$smooth Carmichael-Wintner coefficients, relative to the two classes of functions: $\IPP$ \& $\NSL$. 
\par
We start with the more \lq \lq {\it arithmetic}\rq \rq, so to speak, class, namely the $\IPP$ functions. 
\smallskip
\par
\noindent {\bf Property 1}. {\it Let $F$ be } $\IPP$. {\it Then, } {\stampatello fix a prime} $P$, {\it getting}: 
\smallskip
\item{$(0)$} $\exists \WintnerTP F$, $\exists \CarmichaelTP F$ {\stampatello and } $\CarmichaelTP F=\WintnerTP F$ {\stampatello is square-free supported}
\smallskip
\item{$(1)$} $(\ast)_{(P)}$ : \qquad \qquad \qquad $\forall a\in \N$, \quad ${\displaystyle F_{(P)}(a)=\sum_{{q\in (P)}\atop {q\le \prod_{p\le P}p}}\left( \WintnerPqF \right)c_q(a) }$, 
\par
\noindent
{\it that is} {\stampatello the $P-$local Ramanujan-Wintner expansion, has bounded length} ({\it not $a-$depending})
\smallskip
\item{$(2)$} {\stampatello coefficients in $(\ast)_{(P)}$ are unique }:
\vskip-0.15cm
$$
\exists G_P:\N \rightarrow \C \enspace \hbox{\stampatello with} \enspace \supporto(G_P)\subseteq (P), G_P = \mu^2 \cdot G_P \enspace \hbox{\stampatello and}
$$
\vskip-0.45cm
$$
\forall a\in \N, \quad F_{(P)}(a)=\sum_{{q\in (P)}\atop {q\le \prod_{p\le P}p}}G_P(q)c_q(a)
\leqno{(\ast)\,:\,}
$$
\vskip-0.45cm
\par
\noindent
{\stampatello entail } $G_P=\WintnerTP F$ . 

\medskip

\par
We come to the more \lq \lq {\it analytic}\rq \rq, so to speak, class, namely the $\NSL$ functions. 
\smallskip
\par
\noindent {\bf Property 2}. {\it Let $F$ be } $\NSL$. {\it Then, } {\stampatello fix a prime} $P$, {\it getting}: 
\smallskip
\item{$(0)$} $\exists \WintnerTP F$, $\exists \CarmichaelTP F$, $\CarmichaelTP F=\WintnerTP F$ {\stampatello and for } $q\in (P)$, $\CarmichaelPqF=\WintnerPqF=O_{\varepsilon,P}(q^{-\varepsilon})$
\smallskip
\item{$(1)$} $(\ast)_{(P)}$ : \qquad \qquad \qquad $\forall a\in \N$, \quad ${\displaystyle F_{(P)}(a)=\sum_{q\in (P)}\left( \WintnerPqF \right)c_q(a) }$ 
\par
\noindent
{\it is} {\stampatello the $P-$local Ramanujan-Wintner expansion} 
\smallskip
\item{$(2)$} {\stampatello coefficients in $(\ast)_{(P)}$ are unique }:
\vskip-0.15cm
$$
\exists G_P:\N \rightarrow \C \enspace \hbox{\stampatello with} \enspace \supporto(G_P)\subseteq (P), G_P(q)=O_{\varepsilon,P}(q^{-\varepsilon}), \forall q\in (P) \enspace \hbox{\stampatello and}
$$
\vskip-0.45cm
$$
\forall a\in \N, \quad F_{(P)}(a)=\sum_{q\in (P)}G_P(q)c_q(a)
\leqno{(\ast)\,:\,}
$$
\vskip-0.45cm
\par
\noindent
{\stampatello entail } $G_P=\WintnerTP F$ . 
\smallskip
\par
\noindent
(Hint: both Proofs use that in an absolutely converging double-series we may exchange summations.) 

\vfill

We wish, here, to introduce next section, with new elementary methods.  

\eject

\par				
\noindent{\bf 5. General elementary methods introducing new ideas}
\bigskip
\noindent
We gather some complementary results, having elementary proofs, which supply standard new methods for the study of Ramanujan expansions: especially the ones with smooth summation and/or Wintner coefficients. 
\smallskip
\par
\noindent{\bf 5.1. Ramanujan Clouds}
\smallskip
\par
\noindent
We start with a very easy result that connects absolutely convergent and smooth summation convergent Ramanujan expansions. We recall, for this reason, the notation (compare [C3]) for Ramanujan Clouds : 
$$
\CloudF \defineq \left\{ G:\N \rightarrow \C \left| \right. \forall a\in \N, F(a)=\sum_{q=1}^{\infty}G(q)c_q(a)\right\}
$$
\par
\noindent
is the {\stampatello Ramanujan cloud} of our $F$, namely the set of \lq \lq {\stampatello classic}\rq \rq, say, {\stampatello Ramanujan coefficients}, for a fixed $F:\N \rightarrow \C$; then, we have another set of Ramanujan coefficients for $F$, constituting the {\stampatello Ramanujan smooth cloud} of our $F$ : 
$$
\SmoothCloudF \defineq \left\{ G:\N \rightarrow \C \left| \right. \forall a\in \N, F(a)=\lim_P \sum_{q\in (P)}G(q)c_q(a)\right\}
$$
\par
\noindent
where, in fact, we take (for $P\in \Primes$) the $P-$smooth partial sums' limit over $P\in \Primes$. We complete the notation with the {\stampatello Ramanujan absolute cloud} of our $F$ : 
$$
\AbsCloudF \defineq \left\{ G\in \CloudF \left| \right. \forall a\in \N, \sum_{q\in \N}|G(q)c_q(a)|<\infty\right\}, 
$$
\par
\noindent
the set of classic Ramanujan coefficients of our $F$, in ABSOLUTELY converging $F$ Ramanujan expansions. 
\par
We start noticing that, for the null-function $\0$ we have $<\0> \; \neq \; \subset \0 \supset$: 
$$
G=C\; \hbox{\rm is\enspace constant}
\enspace \Rightarrow \enspace
\sum_{q\in (P)}G(q)c_q(a)=C\prod_{p\le P}\sum_{K=0}^{v_p(a)+1}c_{p^K}(p^{v_p(a)})=\0(a)
\enspace \Rightarrow \enspace
G\in \subset \0 \supset, 
$$
\par
\noindent
compare : Main Lemma in [C2], for the calculation of present $p-$Euler factors (the $K-$sum here). However, a constant function $G\neq \0$ can't be a Ramanujan coefficient of ANY $F:\N \rightarrow \C$, as, for example at $a=1$, we don't have convergence for the \lq \lq classic\rq \rq, say, series : 
$$
C\neq 0, a=1
\enspace \Rightarrow \enspace
\sum_{q=1}^{\infty}G(q)c_q(a)=C\sum_{q=1}^{\infty}\mu(q)
\enspace \hbox{\rm doesn't\enspace converge\enspace in} \enspace \C. 
$$ 
\par
\noindent
(The same coefficients, with summation over $P-$smooth partial sums, give convergence, to $0$ here, see above.) 
\par
In particular, it doesn't converge absolutely, too. We now know that : Ramanujan smooth clouds are NOT contained in Ramanujan absolute clouds (of course for the same $F$). The converse is true since : 
$$
\sum_{q\not\in (P)} |G(q)c_q(a)|\le \sum_{q>P} |G(q)c_q(a)|, 
$$
\par
\noindent
whatever is $G:\N \rightarrow \C$ and $\forall P\in \Primes$. Actually, we have proved that $< \0 >_{\hbox{\stampatello abs}}$ is STRICTLY CONTAINED in $\subset \0 \supset$; for a general $F:\N \rightarrow \C$ it is also true : it follows from the fact that \enspace $\1 \in \SmoothCloud0$ and $G\in \AbsCloudF$ imply $G+\1\in \SmoothCloudF$, but \thinspace $G+\1\not\in \AbsCloudF$. By the way, given any $F$, $\AbsCloudF\not = \emptyset$ because it contains $\hbox{\rm Hil} \thinspace F$, the Hildebrand coefficient ([ScSp], page 166) of our $F$. \enspace In all, we have proved the following. 
\smallskip
\par
\noindent {\bf Proposition 2}. {\it Given any $F:\N \rightarrow \C$, we have $\hbox{\rm Hil} \thinspace F\in \AbsCloudF$} ({\it i.e., $F$ Ramanujan Expansion with Hildebrand Coefficient converges absolutely}) {\it and $\AbsCloudF$ is strictly contained in } $\SmoothCloudF$. 
\smallskip
\par
In particular, we know that {\it all Ramanujan smooth clouds are non-empty}. 

\vfill
\eject

\par				
\noindent{\bf 5.2. Wintner Assumption, Wintner Smooth Assumption and beyond}
\smallskip
\par
\noindent
An even more general hypothesis, starting from $\WA$, than Wintner Smooth Assumption $\WSA$ above, is of course (compare the {\it caveat} soon after $\WSA$ above) the following \lq \lq {\stampatello Wintner Weak Assumption}\rq \rq: 
$$
\exists P_F\in \Primes : \sum_{d\not\in (P)}{{|F'(d)|}\over d}<\infty,\enspace \forall P>P_F. 
\leqno{\WWA}
$$
\par
\noindent
Trivially \enspace $\WA \Rightarrow \WSA \Rightarrow \WWA$. Unexpectedly, for \lq \lq softly decaying\rq \rq, say, Wintner coefficients (compare next $\DD$ in next result with general definition [C3]), we have $\WWA \Rightarrow \WA$. 
\smallskip
\par
\noindent {\bf Proposition 3}. {\it Let $F:\N \rightarrow \C$ have } $\WintnerT F$, {\it with the following, say, } \lq \lq {\stampatello Delange Dual Hypothesis}\rq \rq:
$$
\sum_{q=1}^{\infty}2^{\omega(q)}\left| \Wintner_q F\right|<\infty.
\leqno{\DD}
$$
\par
\noindent
{\it Then, } $\WWA \Rightarrow \WA$. 
\smallskip
\par
\noindent {\bf Proof}. Use, $\forall P\in \Primes$ {\stampatello fixed, the P-orthogonal Wintner Decomposition for } $F'$, i.e. $(7)$ above: 
$$
{{F'(d)}\over d}=\sum_{K\in (P)}\mu(K)\Wintner_{dK} F - \sum_{{r\in )P(}\atop {r>1}}{{F'(dr)}\over {dr}}
\Rightarrow 
\sum_{d\in (P)}{{|F'(d)|}\over d}\le \sum_{d\in (P)}\sum_{K\in (P)}\mu^2(K)\left| \Wintner_{dK} F\right|+\sum_{m\not\in (P)}{{|F'(m)|}\over m}
$$
\par
\noindent
and, passing to \enspace ${\displaystyle \lim_P }$\enspace and applying $\DD$ above with $dK=q$, we get $\WWA \Rightarrow \WA$.\qed

\medskip

\par
On the same lines of Corollary 1, it follows next stronger Corollary: simply from $\BH$ implying finiteness of $\WintnerT\,F$ support (giving $\DD$ trivially). 
\smallskip
\par
\noindent $\WWA-${\bf Corollary 1}. {\it Correlations $F(a):=C_{f,g}(N,a)$ with $\BH$ and $\WWA$ have the $\REEF$.} 

\medskip

\par
Of course, for all fixed $P\in \Primes$, any function $F:\N \rightarrow \C$ satisfying $\IPP$ has 
$$
\sum_{d\in (P)}{{|F'(d)|}\over d}<\infty, 
\enspace \hbox{\rm being\enspace a\enspace finite\enspace sum}, 
$$ 
\par
\noindent
whence it has $\WA$ IFF (if \& only if) it has $\WWA$. The same property, for all fixed $P\in \Primes$, is shared by any $F$ with $\NSL$, from $(1)$. 

\medskip

\par
In view of this last property, since $\BH-$correlations satisfy Ramanujan Conjecture (see [C1] Proposition 1 for this, quoted in $\S1.2$), whence $\NSL$, previous Corollary is actually prefectly equivalent to above Corollary 1. In other words, the difference in between $\WA$ and $\WWA$ may be appreciated only in very general so-to-speak environments for $F$. 

\bigskip

\par
\noindent{\bf 5.3. The {\stampatello Reef} in general}
\smallskip 
\par
\noindent
We saw the applications to $\BH-${\stampatello correlations} (esp., Corollary 1) in $\S1.2$, of our results for general $F$, in order to get the $\REEF$. 
We warn the reader that, in this subsection, $F\neq \0$. Also, see the following, the case $\WintnerT \,F=\0$ is, say, a \lq \lq singular one\rq \rq. 
\smallskip 
\par
We wish to generalize the {\stampatello concept of $\REEF$}, that {\stampatello regards correlations}; for any general $F:\N \rightarrow \C$ we say (notice the notational difference : no dots) 
$$
F \thinspace \hbox{\it has\thinspace the } \Reef
\enspace \definiz \enspace 
\forall a\in \N, \; F(a)=\sum_{q\le Q}\left( \Wintner_q\, F\right)c_q(a), 
$$ 
\par
\noindent
for some FIXED CONSTANT $Q\in \N$. From this property, we get that
$$
F' \thinspace \hbox{\it has\thinspace the } \Reef
\enspace \definiz \enspace 
\forall d\in \N, \; F'(d)=d\sum_{K\le {Q\over d}}\mu(K)\Wintner_{dK}\, F, 
$$ 
\par				
\noindent
for the SAME $Q$ AS ABOVE. In fact, {\stampatello applying Eratosthenes Transform to the}, say, $F-\Reef$, we get the $F'-\Reef$, simply by {\stampatello Kluyver's formula} (after $(2)$ above): 
$$
F(a)=\sum_{d|a}F'(d)=\sum_{q\le Q}\left( \Wintner_q\, F\right)\sum_{{d|a}\atop {d|q}}d\mu(q/d)=\sum_{d|a}d\sum_{K\le {Q\over d}}\mu(K)\Wintner_{dK}\, F, 
$$ 
\par
\noindent
after M\"{o}bius Inversion [T]. 
\par
On the other hand, summing over the divisors $d\in \N$ of $a\in \N$, we get the $F-\Reef$, from the $F'-\Reef$. This idea, of connecting the Ramanujan expansion of a fixed $F$ to an expansion for its Eratosthenes Transform $F'$, goes back to Lucht (see [C3], Proposition 2). 
\par
In particular, the $F-\Reef$ implies that $\WintnerT \,F$ has {\stampatello support} $\supporto(\WintnerT)\subseteq [1,Q]$, apart from the trivial, implicit property: $\exists \WintnerT \,F$. 
\par
On the converse, {\stampatello we ask : once we know that} $\supporto(\WintnerT)\subseteq [1,Q]$, for some $Q\in \N$, {\stampatello under which conditions we get the} $F-\Reef$ {\stampatello above}?
\par
For example, Theorem 1 ensures that $\WA$ gives {\stampatello the Ramanujan-Wintner Smooth expansion}; once we join to this: $\supporto(\WintnerT)$ {\stampatello is finite}, we get the $\Reef$. We similarly prove the following. 
\smallskip
\par
\noindent {\bf Theorem 3}. ({\stampatello $P-$infinitesimal irregular series \& definitively vanishing $\WintnerT$ imply the $\Reef$})
\par
\noindent
{\it Let } $F:\N \rightarrow \C$ {\it have } $\WintnerT\, F$ {\it and } {\stampatello assume } $\IrrPF\to \0$, {\stampatello as } $P\to \infty$ {\stampatello in primes}. {\it Then,}
$$
\left| \supporto(\WintnerT\, F)\right|<\infty
\enspace \Longrightarrow \enspace 
F\enspace \hbox{\it has\enspace the}\enspace \Reef. 
$$
\smallskip
\par
\noindent
We supply a complete and explicit Proof, gathering above properties. 
(Alternatively use Th.2 \& Lemma 3.)
\smallskip
\par
\noindent {\bf Proof}. From Lemma 2, the existence of $\WintnerT\, F$ implies the existence of our $F$ irregular series, compare Remark 3 \& Remark 5. From Lemma 3, the vanishing hypothesis for the irregular series entails (being equivalent to) $\underline{\hbox{\rm the}}$ Ramanujan-Wintner smooth expansion, for $F$; which, under the finiteness for $\WintnerT\, F$ support, implies the $F-\Reef$.\hfil \qed 

\medskip

\par
Of course, the main hypothesis in this Theorem, like also in applications to correlations, is the one for the vanishing of our $F$ irregular series over $P$, as $P\to \infty$ (in the primes). See that, while in previous approaches we rely on less general hypotheses, here a kind of top-generality-hypothesis, so to speak, like this irregular-series-vanishing stops any quest for $\WA$ generalizations, that we briefly described in previous subsection. In fact, $\IrrPF \to \0$ as $P\to \infty$ is EQUIVALENT to THE RWS expansion; this last ingredient only needs finiteness of non-zero Wintner coefficients to produce the $F-\Reef$ as Theorem 3 illustrates. Recall that our previous study, regarding finite Ramanujan expansions [CM], proves that in case of FIXED LENGTH Ramanujan Expansions (like the $\REEF$ \& the $\Reef$ for general $F$) our arithmetic function is a TRUNCATED DIVISOR SUM (with divisors $d\le Q$, for Reefs over $q\le Q$ : see [C3] Theorem 3).  
\smallskip 
\par
For the fixed length Ramanujan expansion $F(a)=\sum_{q}G_F(q)c_q(a)$, $\forall a\in \N$, with coefficients $G_F$, we set 
$$
\ell_F \defineq \sup\{q\in \N : G_F(q)\neq 0\}, 
\qquad \hbox{\rm hereafter\enspace assuming}\enspace G_F\neq \0, 
$$
\par
\noindent
which, of course, is finite IFF the Ramanujan expansion of our $F$ with coefficients $G_F$ has a fixed length; however, it's $+\infty$ IFF such Ramanujan expansion has NOT fixed length. See that, for example, we might have a length depending on $a\in \N$, say $\ell_F(a)$, getting $\ell_F \defineq \sup_{a\in \N}\ell_F(a)$. Compare Theorem 1'[C1] in $\S4$. See that, of course, $\ell_F \in \N$ IFF our F has the $\Reef$, from: $F'$ has the $\Reef$ $\Rightarrow $ $G_F=\WintnerT F$. 
\par
Analogously, for a general $F\neq \0$, we may define \enspace ${\displaystyle d_F\defineq \sup\{d\in \N : F'(d)\neq 0\} }$ \enspace and this, say, \lq \lq top divisor\rq \rq, in finite case (otherwise it's $+\infty$, \lq \lq for almost all arithmetic functions\rq \rq), is linked to $\ell_F$ as $\ell_F=d_F$ (true even in not finite case, as $\ell_F=+\infty=d_F$, then); this follows from quoted Theorem 3 [C3]. See that, of course, $d_F$ is finite IFF our $F$ is a TRUNCATED DIVISOR SUM, with top divisor $d_F\in \N$. Compare next subsection's definition of $Q_F$ in case $F=\0$ : accordingly, we may define $\ell_{\0}\defineq 0$ (but NOT $d_{\0}\defineq 0$ !). 

\medskip

We conclude this brief ride on the $F$-$\Reef$s highlighting the ABSOLUTE CONVERGENCE OF fixed length Ramanujan expansions, whence of THE $F$-$\Reef$. 

\medskip

\par
However, we saw above, there are constant functions $G$ in the Ramanujan smooth cloud of $\0$, while (apart from $G=\0$ itself) there are none in the Ramanujan clouds!
\par
Needless to say, the Panorama of Ramanujan Clouds is very different from Smooth Ramanujan Clouds Landscape$\ldots$ ! 

\vfill
\eject

\par				
\noindent{\bf 5.4. Analytic part and irregular part of arithmetic functions} 
\smallskip 
\par
\noindent
In this subsection, we further generalize previous approach and we {\stampatello study the set}
$$
\C^{\N}_{\hbox{\stampatello fin-win} }\defineq \{ F:\N \rightarrow \C\; |\; \exists \WintnerT \, F \enspace \& \enspace \exists Q\in \N: \supporto(\WintnerT \, F)\subseteq [1,Q]\,\}
$$ 
\par
\noindent
of arithmetic functions {\stampatello with finitely-supported Wintner Transform} (i.e., only a finite number of Wintner coefficients doesn't vanish). For all such $F$ with $\WintnerT \,F\neq \0$, we define  {\stampatello Wintner's range} $Q_F\defineq \sup \supporto(\WintnerT F)$, but this definition also says $Q_F=+\infty$ IFF our $F\not \in \C^{\N}_{\hbox{\stampatello fin-win} }$. While, in case $F$ has $\WintnerT \,F=\0$ we set \enspace $Q_F\defineq 0$. 
In other words, $\forall F\in \C^{\N}_{\hbox{\stampatello fin-win} }$, $Q_F$ is the maximum $q$ with $\Wintner_q\,F\neq 0$.
\smallskip 
\par
Then for these functions $F$, from $\exists \WintnerT\, F$ and $\supporto(\WintnerT\, F)\subseteq [1,Q_F]$, Lemma 2 equation $(8)$ entails
$$
F(a)=\sum_{q\le Q_F}\left(\Wintner_q\, F\right)c_q(a)-\sum_{d|a}\IrrQFdF, 
\quad \forall a\in \N, 
\leqno{\FAI}
$$
\par
\noindent
where now $Q_F\in \N$ {\stampatello might be non-prime}; in this case, we may substitute $Q_F$ with biggest prime $P\le Q_F$, say $P_F$, using the property of the {\stampatello irregular series}, compare Remark 6, of being constant w.r.t. $P\in \Primes$ as long as $P\ge P_F$. \enspace Notice : if $\WintnerT \,F=\0$, then $Q_F=0$ gives the expected empty sum over $q$ inside $\FAI$.

We call this {\stampatello equation} $\FAI$ from the $F=A_F-I_F$ {\stampatello analytic-irregular decomposition of our fixed } $F\in \C^{\N}_{\hbox{\stampatello fin-win} }$, where
$$
A_F(a)\defineq \sum_{q\le Q_F}\left(\Wintner_q\, F\right)c_q(a)
 = \sum_{q\le Q_F}\left(\Wintner_q\, F\right)\sum_{j\in \Z_q^*}e^{2\pi ija/q},
\quad \forall a\in \C
$$
\par
\noindent
is the, {\stampatello say, $F-$analytic part}, that's in fact a Holomorphic function of \enspace $a\in \C$ : $A_F\in \H(\C)$; while, $I_F$ is the, {\stampatello say, $F-$irregular part}, defined $\forall a\in \N$ in terms of {\stampatello irregular series over} the prime $P_F\defineq \max\{ p\in \Primes : p\le Q_F \}$ (we saw above) and we also write $Q_F$ instead of $P_F$ by abuse of notation: 
$$
I_F(a)\defineq \sum_{d|a}\IrrQPdF = \sum_{d|a}\IrrQFdF, 
\quad \forall a\in \N.
$$
\par
\noindent
Inside our {\stampatello fin-win} set of arithmetic functions $F$, previous Theorem 3 is now very clear: the $F-\Reef$ is equivalent to having $I_F=\0$ in $\FAI$! This also reveals that the functions $F$ in our set, having the $\Reef$, are entire functions and, by Liouville Theorem, $F(a)$ is bounded $\forall a\in \C$ IFF our $F$ is a constant !! So, once again (compare quoted property from [CM]) a kind of \lq \lq rarity\rq \rq, say, is the $F-\Reef$ !!!  

\medskip

\par
Thanks to $\FAI$ we might think about $\underline{\rm the}$ Ramanujan-Wintner Smooth expansion, say RWSE, for a fixed ARBITRARY $F\in \C^{\N}$, as a process of asymptotic approximations, as $Q\to \infty$, by functions $F\in \C^{\N}_{\hbox{\stampatello fin-win} }$, each with Wintner coefficients vanishing after $Q_F$ ! From this point of view, $\FAI$, itself, is a kind of \lq \lq {\stampatello approximate Reef}\rq \rq. (Compare page 8 in [C1, version 3].)

\medskip

\par
See that having $Q_F$ doesn't suffice to get the $\Reef$ for $F$ (compare Counterexample 1 [C1] in $\S5.6$). We saw in previous subsection that: $F$ has the $\Reef$ $\Leftrightarrow $ $\ell_F=d_F$ are {\stampatello both finite}. Is there a condition under which we get the $\Reef$ for $F$, when $Q_F\in \N$ ? Next result answers. 
\smallskip
\par
\noindent {\bf Theorem 4}. {\it Let } $F:\N \rightarrow \C$ {\it have finite } $Q_F$. {\it Then}
$$
F\enspace \hbox{\stampatello has\enspace the\enspace RWSE}
\enspace
\Longleftrightarrow 
\enspace
F\enspace \hbox{\stampatello has\enspace the\enspace }\Reef
\enspace
\Longleftrightarrow 
\enspace
F\enspace \hbox{\stampatello satisifies\enspace the\enspace}\WA. 
$$
\smallskip
\par
\noindent
We leave the Proof to the interested reader. 
\smallskip
\par
\noindent {\bf Remark 11.} In case $Q_F=+\infty$, we may even have $F$ {\it with} the RWSE, but {\it without} the $\WA$. (Hint: esp., $F_{\alpha}$ having a completely multiplicative $F_{\alpha}'$, with $F_{\alpha}'(p):=e(\alpha p)$, $\forall p\in \Primes$, with a fixed irrational $0<\alpha<1$).\hfill $\diamond$ 

\medskip

Last, but not least, notice that, when coming back to correlations, the Basic Hypothesis makes our $F(a):=C_{f,g}(N,a)$ (again, from Proposition 1 [C1,version 3]) have a finitely-supported Wintner Transform, i.e. : $F\in \C^{\N}_{\hbox{\stampatello fin-win} }$. Thus $\FAI$ can turn into a practical \& effective formula for estimating the remainder for $-I_F=F-A_F$, where $A_F$ now's nothing but the fixed-length Ramanujan expansion in the $\REEF$ ! In other words, even if we might, as it's \lq \lq too rare\rq \rq, not have the $\REEF$, we might, from $\FAI$, try to estimate the remainder, in terms of our $F-${\stampatello irregular series} (in $I_F$, here), for the \lq \lq {\stampatello Hardy-Littlewood asymptotics}\rq \rq, for general $\BH-$correlations $F$. (Compare [C0] formul\ae.) 

\vfill
\eject

\par 				
\noindent{\bf 5.5. Irregular series of multiplicative arithmetic functions} 
\smallskip 
\par
\noindent
The Irregular Series, $\IrrPdF$, for general $F$ having $\WintnerT\,F$, even assuming all the hypotheses above, remains a kind of mistery. 
\par
However, for $F:\N \rightarrow \C$ a {\stampatello multiplicative} arithmetic function, it simplifies a lot, as we see now: 
$$
F\enspace \hbox{{\stampatello multiplicative}} 
\enspace \Rightarrow \enspace 
\IrrPdF=F'(d)\enspace \IrrP1F, 
\quad \forall P\in \Primes, \forall d\in (P). 
$$

\bigskip
\bigskip
\bigskip

\par
\noindent{\bf 5.6. Correlations with Basic Hypothesis, but without Reef: studying Counterexample 1} 
\smallskip 
\par
\noindent
The Counterexample 1 [C1] shows that the Basic Hypothesis, implying that the Wintner transform is finitely-supported, is NOT sufficient to get the Reef. 
\par
We recall briefly that Counterexample 1, see [C1], is the correlation of two arithmetic functions, say  $f_0,g_0:\N \rightarrow \C$, chosen this way: 
\par
\noindent 
FIX $N,Q\in \N$ with $Q\le N$ and two integers $1\le n_0\le N$ and $2<q_0\le Q$. Choose 
$$
f_0(n)\defineq \1_{\{ n_0\}}(n), 
\enspace \forall n\in \N
\quad 
\hbox{\rm and}
\quad
g_0(m)\defineq c_{q_0}(m), 
\enspace \forall m\in \N 
$$
\par
\noindent 
whence : 
$$
n_0\equiv -1(\bmod q_0) 
$$
\par
\noindent 
implies : we can't have the Reef for $C_{f_0,g_0}(N,a)$, $\forall a\in \N$, since in particular for $a=1$ Reef's LHS and RHS are DIFFERENT. (See page 8,[C1], for details). 
\par
We profit, here, to gather some properties of our Counterexample 1, we'll call $F_0(a)$, in the more manageable case that modulus $q_0$ is prime $q_0=p_0\in \Primes$ and $n_0=q_0-1=p_0-1$: hereafter, with \enspace $p_0>2$, 
$$
F_0(a)\defineq c_{p_0}(a-1),
\qquad
\forall a\in \N. 
$$
\par
\noindent 
(It might seem that this is not a correlation, but please gather above definitions!) 
\par
This correlation satisfies, as usual, our Basic Hypothesis and, by the way, has Wintner Transform $\WintnerT F_0$ simply given by $q-$th coefficient ${1\over {\varphi(p_0)}}c_{p_0}(p_0-1)={1\over {\varphi(p_0)}}\mu(p_0)=-1/\varphi(p_0)$, if and only if $q=p_0$, vanishing otherwise. (In particular, this Tranform is finitely-supported, of course.) We start calculating Eratosthenes Transform: 
$$
F_0'(1)=F_0(1)=\varphi(p_0)=p_0-1,
\qquad
\hbox{\rm while}
\qquad
d>1 \enspace \Rightarrow \enspace 
F_0'(d)=p_0 S_{p_0}(d), 
$$
\par
\noindent 
where we set
$$
S_{p_0}(d)\defineq \sum_{{a|d}\atop {a\equiv 1\bmod p_0}}\mu\left({d\over a}\right),
\qquad
\forall d\in \N, 
$$
\par
\noindent 
because : $\forall d\in \N$ we have 
$$
F_0'(d)=\sum_{{a|d}\atop {a\equiv 1\bmod p_0}}\varphi(p_0)\mu\left({d\over a}\right)+\sum_{{a|d}\atop {a\not \equiv 1\bmod p_0}}\mu(p_0)\mu\left({d\over a}\right)
=(p_0-1)\sum_{{a|d}\atop {a\equiv 1\bmod p_0}}\mu\left({d\over a}\right)-\sum_{{a|d}\atop {a\not \equiv 1\bmod p_0}}\mu\left({d\over a}\right),  
$$
\par
\noindent 
which is \enspace $p_0 S_{p_0}(d)$, $\forall d>1,$ from M\"{o}bius inversion: 
$$
\sum_{{a|d}\atop {a\not \equiv 1\bmod p_0}}\mu\left({d\over a}\right)=\sum_{a|d}\mu\left({d\over a}\right)-\sum_{{a|d}\atop {a\equiv 1\bmod p_0}}\mu\left({d\over a}\right)
 =-\sum_{{a|d}\atop {a\equiv 1\bmod p_0}}\mu\left({d\over a}\right), 
\enspace \forall d>1. 
$$ 
\par
\noindent 
Dirichlet characters modulo $p_0$ allow to write
$$
S_{p_0}(d)={1\over {\varphi(p_0)}}\sum_{\chi(\!\!\bmod p_0)}\sum_{a|d}\chi(a)\mu\left({d\over a}\right),
\qquad
\forall d\in \N\backslash \{1\}. 
\leqno{(\ast)_{p_0}}
$$

\vfill
\eject

\par 				
\noindent
We may distinguish three cases, for the integers $d>1$ : 
$$
v_{p_0}(d)=0
\quad \hbox{\rm and} \quad 
d>1
\leqno{({\bf 0})}
$$
$$
v_{p_0}(d)=1
\leqno{({\bf 1})}
$$
$$
v_{p_0}(d)\ge 2
\leqno{({\bf 2})}
$$
\par
\noindent 
Last case $({\bf 2})$ is the simplest, since, setting $K:=d/p^{v_{p_0}(d)}\in \Z_{p_0}^*$, 
$$
S_{p_0}(d)=S_{p_0}(p_0^{v_{p_0}(d)}\cdot K)
 =\sum_{{a|K}\atop {a\equiv 1\bmod p_0}}\mu\left(p_0^{v_{p_0}(d)}\cdot {K\over a}\right)
  =\mu\left(p_0^{v_{p_0}(d)}\right)S_{p_0}(K)
   =0, 
$$
\par
\noindent 
in case $({\bf 2})$. 
\par
Similarly, setting in case $({\bf 1})$ $K:=d/p_0\in \Z_{p_0}^*$, 
$$
S_{p_0}(d)=S_{p_0}(p_0\cdot K)
 =\sum_{{a|K}\atop {a\equiv 1\bmod p_0}}\mu\left(p_0\cdot {K\over a}\right)
  =\mu(p_0)S_{p_0}(K)
   =-S_{p_0}(d/p_0), 
$$
\par
\noindent 
in case $({\bf 1})$. In particular, we may omit the single $d=p_0$, as $S_{p_0}(p_0)=-S_{p_0}(1)=-1$. 
\par
Everything boils down to case $({\bf 0})$, in which formula $(\ast)_{p_0}$ at previous page, with Dirichlet characters, becomes: 
$$
S_{p_0}(d)={1\over {\varphi(p_0)}}\sum_{\chi(\!\!\bmod p_0)}\sum_{K|d}\mu(K)\chi(d)\overline{\chi}(K)
 ={1\over {\varphi(p_0)}}\sum_{\chi(\!\!\bmod p_0)}\chi(d)\prod_{p|d}\left(1-\overline{\chi}(p)\right), 
$$
\par
\noindent 
in case $({\bf 0})$, because the flipping $K:={d\over a}$ of divisors $a|d$ has 
$$
\chi\left({d\over K}\right)={{\chi(d)}\over {\chi(K)}}
 =\chi(d)\overline{\chi}(K),
\quad
\forall K|d \enspace (\hbox{\rm recall}\enspace d\in \Z_{p_0}^*).
$$
\par
\noindent 
In this formula, the finite product over primes $p$ dividing $d$ (from : $p\equiv 1\bmod p_0$ $\Rightarrow $ $\overline{\chi}(p)=1$, $\forall \chi\bmod p_0$) immediately entails the property 
$$
\exists p|d \enspace : \enspace p\equiv 1\bmod p_0
\quad \Rightarrow \quad
S_{p_0}(d)=0. 
\leqno{(\ast)_0}
$$
\par
\noindent 
We may so to speak summarize these properties of $F_0'$, giving a glance to (without calculating it) the mean value of $|F_0'|$, i.e.: 
$$
\lim_x {1\over x}\sum_{d\le x}|F_0'(d)|=p_0\lim_x {1\over x}\left(\sum_{{1<d\le x}\atop {(d,p_0)=1}}|S_{p_0}(d)|+\sum_{{1<d\le x/p_0}\atop {(d,p_0)=1}}|S_{p_0}(d)|\right)
=
$$
$$
=p_0\lim_x {1\over x}\left(\sumflat_{1<d\le x}|S_{p_0}(d)|+\sumflat_{1<d\le {x\over {p_0}}}|S_{p_0}(d)|\right)
 =(p_0+1)\lim_x {1\over x}\sumflat_{1<d\le x}|S_{p_0}(d)|, 
$$
\par
\noindent 
where the first equation comes from distinguishing cases $(\bf 0)$ and $(\bf 1)$, while second one introduces the $\flat$ notation, from $(\ast)_0$ property, that means: any prime $p|d$ is NEITHER $0$ NOR $1$ modulo $p_0$; finally, last equation, so to speak, comes from the change of variable in second limit passing from $x$ to $p_0x$. 

\vfill
\eject

\par				
\noindent
We still have two properties of our $F_0$ that are noteworthy to see: namely, we give a brief look at the behavior of $S_{p_0}(d)$, respectively on square-free $d>1$ and on the powers of primes $p$ different from $p_0$. 
\par
First of all, see that on square-free $d>1$ we have
$$
S_{p_0}(d)=\mu(d)\sum_{{a|d}\atop {a\equiv 1\bmod p_0}}\mu(a), 
$$
\par
\noindent 
from the trivial remark that these $d$ have $\mu(d/a)=\mu(d)\mu(a)$, because $a$ is square-free, too, and $1/\mu(a)=\mu(a)$ in this case. Hence,
$$
\left|S_{p_0}(d)\right|=\left|\sum_{{a|d}\atop {a\equiv 1\bmod p_0}}\mu(a)\right|, 
\quad \forall d>1, \enspace \mu^2(d)=1. 
$$
\par
\noindent 
This may be of some help in above calculations for $\left|S_{p_0}(d)\right|$ averages; also, Dirichlet characters modulo $p_0$ simplify above $(\ast)_{p_0}$ as: 
$$
{\widetilde{S}}_{p_0}(d)\defineq \sum_{{a|d}\atop {a\equiv 1\bmod p_0}}\mu(a)
 ={1\over {\varphi(p_0)}}\sum_{\chi(\!\!\bmod p_0)}\prod_{p|d}\left(1-\chi(p)\right),
\enspace 
\forall d\in \N, \mu^2(d)=1. 
\leqno{(\widetilde{\ast})_{p_0}}
$$
\par
\noindent 
Above cases $({\bf 0})$, $({\bf 1})$ and $({\bf 2})$ for $S_{p_0}$ become, for ${\widetilde{S}}_{p_0}(d)$ on square-free $d>1$, only the two possibilities
$$
v_{p_0}(d)=0
\quad \hbox{\rm and} \quad 
d>1, \mu^2(d)=1
\leqno{(\widetilde{0})}
$$
$$
v_{p_0}(d)=1
\quad \hbox{\rm and} \quad 
\mu^2(d)=1
\leqno{(\widetilde{1})}
$$
\par
\noindent 
becoming, for ${\widetilde{S}}_{p_0}$, on the same lines as above, in only one occurrence : 
\par
in case $(\widetilde{1})$, setting $K:=d/p_0\in \Z_{p_0}^*$, to get 
$$
{\widetilde{S}}_{p_0}(d)={\widetilde{S}}_{p_0}(p_0\cdot K)
 =\sum_{{a|K}\atop {a\equiv 1\bmod p_0}}\mu(p_0\cdot K)
  =\mu(p_0){\widetilde{S}}_{p_0}(K)
   =-{\widetilde{S}}_{p_0}(d/p_0). 
$$

\bigskip

Turning back to our $F_0$, we prove now that its Eratosthenes Transform $F_0'(d)$ is NOT infinitesimal as $d\to \infty$; simply, calculating $S_{p_0}(d)$ on $d=p^K$, powers, with infinitely many $K\in \N$, of primes $p\neq p_0$ with $p\not \equiv 1(\bmod \enspace p_0)$, it follows, from next formula, that $F_0'(p^K)=\pm p_0$, for infinitely many $K\in \N$, because: 
$$
S_{p_0}(p^K)=\sum_{{j=0}\atop {p^j\equiv 1(\!\!\bmod p_0)}}^{K} \mu(p^{K-j})
 = \1_{p^K\equiv 1(\!\!\bmod p_0)}-\1_{p^{K-1}\equiv 1(\!\!\bmod p_0)}
  \not \to 0, 
\quad 
\hbox{\rm as }\enspace K\to \infty, 
$$
\par
\noindent 
from the definition of $S_{p_0}$ above (recalling $p_0>2$ here), since Fermat's little Theorem implies that it's $1$ on the $K\equiv 0(\!\bmod \enspace  p_0-1)$ and $-1$ on the $K\equiv 1(\!\bmod \enspace p_0-1)$. Of course, these give two subsequences for $F_0'(d)$ not infinitesimal on $d=p^K$, as $d\to \infty$. 
\medskip
\par
In particular, saying that $F_0'(d)$ doesn't go to $0$ as $d\to \infty$ proves once again that the Reef doesn't hold: in fact, the Reef holds if and only if our Eratosthenes Transform has finite support! 
\medskip
\par
We will study in deeper details : in order to prove whether $\ETD$ holds or not for our $F_0$ above (recall, an instance of [C1] Counterexample 1) we found some technical difficulties we hope to overcome in the future. 
\medskip
\par
Last but not least, we propose an exercise to interested readers. Recall $p_0>2$ in the above $F_0$ definition. 
\smallskip
\par
\noindent {\bf Curiosity 1}. {\it Our } $F_0$ {\it is not } $\IPP$: {\it taking } $p\equiv-1(\bmod \, p_0)$, {\it with } $p>p_0$, {\it we have, when } $a=p^2$,
$$
F_0(\kappa(a))=F_0(p)=-1\neq p_0-1=F_0(p^2)=F_0(a). 
$$ 

\vfill
\eject

\par				
\noindent{\bf 5.7. Smooth/classic Carmichael-Wintner coefficients for imaginary exponentials \& applications} 
\smallskip 
\par
\noindent
We make, say, a kind of exercises in computing resp., the Classic Carmichael $\CarmichaelT F_{j,q}$ and all the Smooth Carmichael $\CarmichaelTP F_{j,q}$ coefficients, $\forall P\in \Primes$, for the remarkable $F_{j,q}(a):=e_q(ja)$ : the imaginary exponential function, where the two parameters $q\in \N$ and $j\in \Z_q^*$ are FIXED. We'll use this notation, recalling: from Lemma 6, since our $F_{j,q}$ is bounded (whence, $\NSL$, too), we have $\WintnerTP F_{j,q}=\CarmichaelTP F_{j,q}$ ($\forall P\in \Primes$). We also know, after finding $\CarmichaelT F_{j,q}$, that it equals $\WintnerT F_{j,q}$, from Proposition 3 in [C3], that is a kind of reformulation of a 1987 result of Delange (see [C3] for the bibliography). We'll indicate $F_{j,q}=e_q(j\;\bullet)\in \C^{\N}$. 
\medskip
\par
Thus 
\smallskip
\par
\noindent {\bf Lemma 7}. ({\stampatello Carmichael coefficients of imaginary exponential function})
\par
\noindent
{\it Fix } $q\in \N$ {\it and } $j\in \Z_q^*$. {\it Then,} $\forall \ell \in \N$, 
$$
\Carmichael_{\ell}\thinspace e_q(j\;\bullet)=\1_{\ell=q}\cdot {1\over {\varphi(q)}}. 
$$
\smallskip
\par
\noindent {\bf Proof}. Carmichael coefficient definition and Kluyver formula (see soon after $(2)$ above)
$$
\Carmichael_{\ell}\enspace e_q(j\;\bullet)={1\over {\varphi(\ell)}}\lim_x {1\over x}\sum_{a\le x}e_q(ja)c_{\ell}(a)
 = {1\over {\varphi(\ell)}}\sum_{d|\ell}d\mu\left({{\ell}\over d}\right)\lim_x {1\over x}\sum_{m\le x/d}e_q(jdm), 
$$
\par
\noindent
with the cancellation in exponential sums, i.e., as $x\to \infty$, 
$$
\sum_{m\le x/d}e_q(jdm)=\1_{d\equiv 0(\!\!\bmod q)}\left[{x\over d}\right]+\1_{d\not \equiv 0(\!\!\bmod q)}O\left({1\over {\left\Vert{{jd}\over q}\right\Vert}}\right)
 = \1_{d\equiv 0(\!\!\bmod q)}\cdot {x\over d}+O_q(1), 
$$
\par
\noindent
give soon the thesis 
$$
\Carmichael_{\ell}\enspace e_q(j\;\bullet)={1\over {\varphi(\ell)}}\sum_{{d|\ell}\atop {d\equiv 0(\!\!\bmod q)}}\mu\left({{\ell}\over d}\right)
 ={1\over {\varphi(\ell)}}\,\cdot\,\1_{q|\ell}\,\cdot\, \sum_{d'\left|{{\ell}\over q}\right.}\mu\left({{\ell/q}\over {d'}}\right)
  = \1_{\ell=q}\,\cdot\,{1\over {\varphi(\ell)}}, 
$$
\par
\noindent
by M\"obius inversion (quoted after $(2)$ above).\hfill \square
\medskip
\par
\noindent {\bf Remark 12.} The main idea is the {\bf resonance of moduli} $q$ and $\ell$. Writing $c_{\ell}(a)$ with the exponentials and applying soon exponential sums cancellation, as an alternative proof, renders this more transparent.\hfill $\diamond$ 

\bigskip

This result is so easy that we may have called it a \lq \lq Fact\rq \rq. In case of our $F_{j,q}$ transform $\CarmichaelTP F_{j,q}$, we need more small ideas combined together: the main anthem is a kind of writing averages over $P-$smooth numbers involving imaginary exponentials in term of same averages over Dirichlet characters, that have a multiplicative structure, instead. 
\medskip
\par
In fact, we start calculating $P-$smooth Carmichael coefficients of a general class of arithmetic functions $F$, the $\NSL$ ones, in terms of $P-$smooth numbers averages, with Ramanujan sums; this will be applied to our imaginary exponential function $F=F_{j,q}$, but the following result is quite general. Proof follows [C1]. 
\smallskip
\par
\noindent {\bf Lemma 8}. ({\stampatello Carmichael $P-$smooth coefficients of $\NSL$ functions})
\par
\noindent
{\it Let } $F:\N\rightarrow \C$ {\it be } $\NSL$. {\it Then,} $\forall P \in \Primes$, 
$$
\Carmichael_{\ell}^{(P)}F={1\over {\varphi(\ell) {\displaystyle \sum_{m\in (P)} }{1\over m}}}\cdot \sum_{t\in (P)}{{F(t)}\over t}c_{\ell}(t),
\quad
\forall \ell \in (P). 
$$
\smallskip
\par
\noindent {\bf Proof}. Carmichael $P-$smooth $\ell-$th coefficient definition and Lemma 1 of [C1] (\lq \lq M\"obius Switch\rq \rq) give  
$$
\Carmichael_{\ell}^{(P)}F={1\over {\varphi(\ell)}}\lim_x {1\over x}\sum_{a\le x}c_{\ell}(a)\sum_{{d\in (P)}\atop {d|a}}F'(d)
 = {1\over {\varphi(\ell)}}\lim_x {1\over x}\sum_{a\le x}c_{\ell}(a)\sum_{{t\in (P)}\atop {{t|a} \atop {{a\over t}\in )P(}}}F(t), 
$$
\par				
\noindent
where the sums exchange, the property $\ell\in (P),m\in )P($ $\Rightarrow$ $c_{\ell}(tm)=c_{\ell}(t)$ and Lemma 2 [C1], a kind of Eratosthenes-Legendre sieve, give 
$$
\sum_{a\le x}c_{\ell}(a)\sum_{{t\in (P)}\atop {{t|a} \atop {{a\over t}\in )P(}}}F(t)=\sum_{t\in (P)}F(t)\sum_{{m\le x/t}\atop {m\in )P(}}c_{\ell}(tm)
 =\sum_{t\in (P)}F(t)c_{\ell}(t)\sum_{{m\le x/t}\atop {m\in )P(}}1
$$
$$
=\sum_{t\in (P)}F(t)c_{\ell}(t)\left(\prod_{p\le P}\left(1-{1\over p}\right){x\over t}+O_P(1)\right)
 =\prod_{p\le P}\left(1-{1\over p}\right){x\over t}\sum_{t\in (P)}F(t)c_{\ell}(t)+O_{P,\ell,F}(1)
$$
\par
\noindent
and recalling (for details, see [C1]: Proposition 2 Proof start)
$$
\prod_{p\le P}\left(1-{1\over p}\right)=\left(\prod_{p\le P}\left(1-{1\over p}\right)^{-1}\right)^{-1}
 =\left(\sum_{m\in (P)}{1\over m}\right)^{-1}
  ={1\over { {\displaystyle \sum_{m\in (P)} }{1\over m}}}
$$
\par
\noindent
gives at once 
$$
\lim_x {1\over x}\sum_{a\le x}c_{\ell}(a)\sum_{{t\in (P)}\atop {{t|a} \atop {{a\over t}\in )P(}}}F(t)={1\over { {\displaystyle \sum_{m\in (P)} }{1\over m}}}\cdot \sum_{t\in (P)}{{F(t)}\over t}c_{\ell}(t), 
$$
\par
\noindent
whence the thesis.\hfill \square
\medskip
\par
Next Lemma is a Corollary of previous one, plus a switch of harmonics: from imaginary exponentials to Dirichlet characters. Gauss sums $\tau(\chi)$  definition [D] is recalled in the Proof. 
\smallskip
\par
\noindent {\bf Lemma 9}. ({\stampatello imaginary exponentials' Carmichael $P-$smooth coeff.s: switch to characters})
{\it Fix } $q\in \N$ {\it and } $j\in \Z_q^*$. {\it Then,} $\forall P \in \Primes$, {\it with } $P\ge q$, 
$$
\Carmichael_{\ell}^{(P)}\thinspace e_q(j\;\bullet)={1\over {\varphi(\ell) {\displaystyle \sum_{m\in (P)} }{1\over m}}}\cdot 
\sum_{{b|q}\atop {q':=q/b}}{1\over {b\varphi(q')}}\sum_{\chi(\!\!\bmod q')}\tau(\overline{\chi})\chi(j)\sum_{t\in (P)}{{\chi(t)}\over t}c_{\ell}(bt),
\quad 
\forall \ell \in (P). 
$$
\smallskip
\par
\noindent {\bf Proof}. Straight from previous Lemma for $F=F_{j,q}=e_q(j\;\bullet)$, 
$$
\Carmichael_{\ell}^{(P)}\thinspace e_q(j\;\bullet)={1\over {\varphi(\ell) {\displaystyle \sum_{m\in (P)} }{1\over m}}}\cdot \sum_{t\in (P)}{{e_q(jt)}\over t}c_{\ell}(t),
\quad
\forall \ell \in (P). 
$$
\par
\noindent
We switch from imaginary exponentials to Dirichlet characters of modulus $q':=q/b$ by the inversion formula (see [D]) with the {\it Gauss sum}
$$
\tau(\chi)\defineq \sum_{m\in \Z_{q'}^*}\chi(m)e_{q'}(m)
\quad \Rightarrow \quad
e_{q'}(k)={1\over {\varphi(q')}}\sum_{\chi(\!\bmod q')}\tau(\overline{\chi})\chi(k),
\enspace \forall k\in \Z_{q'}^*  
$$
\par
\noindent
giving at once, from hypothesis $P\ge q$ entailing $b\in (P)$\enspace $\forall b|q$, the following:  
$$
\sum_{t\in (P)}{{e_q(jt)}\over t}c_{\ell}(t)=\sum_{b|q}\sum_{t\in (P),(t,q)=b}{{e_{q/b}(j(t/b))}\over t}c_{\ell}(t)
 =\sum_{b|q}{1\over b}\sum_{{t'\in (P)}\atop {(t',q/b)=1}}{{e_{q/b}(jt')}\over {t'}}c_{\ell}(bt')=
$$
$$
=\sum_{{b|q}\atop {q':=q/b}}{1\over b}\sum_{{t\in (P)}\atop {(t,q')=1}}{{e_{q'}(jt)}\over t}c_{\ell}(bt)
 =\sum_{{b|q}\atop {q':=q/b}}{1\over {b\varphi(q')}}\sum_{\chi(\!\bmod q')}\tau(\overline{\chi})\chi(j)\sum_{t\in (P)}{{\chi(t)c_{\ell}(bt)}\over t}, 
$$
\par
\noindent
from $j\in \Z_{q'}^*$ and the property: $(t,q')=1$ is implicit in presence of $\chi(t)$, whence the formula.\hfill \square

\vfill
\eject

\par 				
\noindent
We have a kind of two small problems to face, for an explicit formula in terms of characters and partial Euler products. First, we have to get rid of the \lq \lq extra factor\rq \rq, so to speak,in the Ramanujan sum of modulus $\ell$ in the above formula: we solve this in next Lemma, with a small idea (we will \lq \lq kill $b$\rq \rq, say). 
\smallskip 
\par
\noindent {\bf Lemma 10}. ({\stampatello absorbing extra factors in Ramanujan sums})
{\it Choose any } $\ell,b,t\in \N$. {\it Then}
$$
c_{\ell}(bt)={{\varphi(\ell)}\over {\varphi(\ell/(\ell,b))}}\enspace c_{\ell/(\ell,b)}(t). 
$$
\smallskip
\par
\noindent {\bf Proof}. Ramanujan sums Explicit Formula [M, page 22 : H\"older's 1936 formula], applied twice: 
$$
c_{\ell}(bt)=\varphi(\ell)\cdot {{\mu(\ell/(\ell,bt))}\over {\varphi(\ell/(\ell,bt))}}
 =\varphi(\ell)\cdot {{\mu(\ell'/(\ell',b't))}\over {\varphi(\ell'/(\ell',b't))}}
  ={{\varphi(\ell)}\over {\varphi(\ell')}}c_{\ell'}(b't)
   ={{\varphi(\ell)}\over {\varphi(\ell')}}c_{\ell'}(t), 
$$
\par
\noindent
using now $\ell/(\ell,bt)={{\ell/(\ell,b)}\over {(\ell/(\ell,b),tb/(\ell,b))}}=\ell'/(\ell',b't)$, where
$$
\ell':=\ell/(\ell,b),\enspace b':=b/(\ell,b), 
$$
\par
\noindent
together with $b'\in \Z_{\ell'}^*$.\hfill \square

\bigskip

Just like we have, say, separated $b$ from other variables, we need now to separate the prime factors of a fixed modulus $q'$ from other variables, in next Lemma with Dirichlet characters modulo $q'$. In fact, when we want to \lq \lq flip\rq \rq, say, a Dirichlet character $\chi(d)$, over divisors $d|n$, into $\chi(n/K)$, with complementary divisor $K:=n/d$, we may then write $\chi(n/K)=\chi(n)/\chi(K)$ only if we know that $K$ is coprime to $q'$ (our $\chi$ modulus); in other words, we have to separate the prime-factors of $n$ dividing modulus $q'$. As we see soon. 
\smallskip
\par
\noindent {\bf Lemma 11}. ({\stampatello separating modulus prime-factors before flipping Dirichlet characters})
\par
\noindent
{\it Choose any } $\ell',q'\in \N$. {\it Then, setting } $q'':=\prod_{p|\ell',p|q'}p^{v_p(\ell')}$, $\ell'':=\ell'/q''$, {\it we have } $\forall \chi(\!\bmod\enspace q')$ 
$$
\chi'(\ell')=\sum_{d|\ell'}\chi(d)\mu\left({{\ell'}\over d}\right)=\mu(q'')\chi(\ell'')\prod_{p|\ell''}(1-\overline{\chi}(p)). 
$$
\smallskip
\par
\noindent {\bf Proof}. In the sum over $d$, in LHS, the factor $\chi(d)$ implies $(d,q')=1$ and $\ell''\in \Z_{q'}^*$ by construction : 
$$
\sum_{d|\ell'}\chi(d)\mu\left({{\ell'}\over d}\right)=\sum_{{d|\ell'}\atop {(d,q')=1}}\chi(d)\mu\left({{\ell'}\over d}\right)
 =\sum_{d|\ell''}\chi(d)\mu\left(q''\cdot {{\ell''}\over d}\right)
  =\mu(q'')\sum_{K|\ell''}\mu(K)\chi\left({{\ell''}\over K}\right)=
$$
$$
=\mu(q'')\chi(\ell'')\sum_{K|\ell''}\mu(K)\overline{\chi}(K)
 =\mu(q'')\chi(\ell'')\prod_{p|\ell''}(1-\overline{\chi}(p)), 
$$
\par
\noindent
flipping, say, the divisors $d$ as:\enspace $K:=\ell''/d$, having used $K|\ell''$ $\Rightarrow$ $\chi(K)\neq 0$ $\Rightarrow$ $\chi(\ell''/K)=\chi(\ell'')\overline{\chi(K)}$ and the general formula [T] 
$$
\sum_{d|n}\mu(d)f(d)=\prod_{p|n}(1-f(p)),
$$
\par
\noindent
for all multiplicative functions $f$.\hfill \square

\vfill
\eject

\par 				
\noindent
Before gathering all these Lemmas together to compute Carmichael $P-$smooth coefficients of our imaginary exponential function, in next Theorem, we need to look at the corresponding Carmichael coefficients: we express them as the $\chi=\chi_0$ part of Lemma 9 formula, for ALL the principal characters modulo $q'$, $\forall q'\in \N$ (they're the only $\chi$ modulo $q'$, of course, in cases $q'=1,2$). 
\smallskip
\par
\noindent {\bf Lemma 12}. ({\stampatello imaginary exponentials' Carmichael coefficients: principal characters})
\par
\noindent
{\it Fix } $q\in \N$ {\it and } $j\in \Z_q^*$. {\it Then,} $\forall P \in \Primes$, {\it with } $P\ge q$, 
$$
\Carmichael_{\ell}\thinspace e_q(j\;\bullet)={1\over {\varphi(\ell) {\displaystyle \sum_{m\in (P)} }{1\over m}}}\cdot 
\sum_{{b|q}\atop {q':=q/b}}{1\over {b\varphi(q')}}\tau(\overline{\chi_0})\chi_0(j)\sum_{t\in (P)}{{\chi_0(t)}\over t}c_{\ell}(bt),
\quad 
\forall \ell \in (P). 
$$
\smallskip
\par
\noindent {\bf Proof}. Straight from: $\tau(\overline{\chi_0})=\tau(\chi_0)=c_{q'}(1)=\mu(q')$ [D], $\mu(q')/\varphi(q')=c_q(bt)/\varphi(q)$ (from quoted H\"older 1936 formula, in [M]) and $\chi_0(j)=1$ (recall $(j,q)=1=(j,q')$ here), rendering RHS 
$$
{1\over {\varphi(\ell) {\displaystyle \sum_{m\in (P)} }{1\over m}}}\cdot 
\sum_{{b|q}\atop {q':=q/b}}{{c_q(bt)}\over {b\varphi(q)}}\sum_{t\in (P)}{{\chi_0(t)}\over t}c_{\ell}(bt)={1\over {\varphi(\ell) {\displaystyle \sum_{m\in (P)} }{1\over m}}}\cdot 
\sum_{b|q}{{c_q(bt)}\over {\varphi(q)}}\sum_{{t\in (P)}\atop {(bt,q)=b}}{{c_{\ell}(bt)}\over {bt}}, 
$$
\par
\noindent
because : $\chi_0(t)\neq 0$ $\Leftrightarrow $ $(t,q/b)=1$ $\Leftrightarrow $ $(bt,q)=b$, while our RHS is, from $P\ge q$ $\Rightarrow $ $q\in (P)$ $\Rightarrow $ $b\in (P)$, $\forall b|q$, 
$$
{1\over {\varphi(\ell) {\displaystyle \sum_{m\in (P)} }{1\over m}}}\cdot 
{1\over {\varphi(q)}}\sum_{b|q}\sum_{{t\in (P)}\atop {(bt,q)=b}}{{c_{\ell}(bt)c_q(bt)}\over {bt}}
={1\over {\varphi(q)}}\cdot {1\over {\varphi(\ell) {\displaystyle \sum_{m\in (P)} }{1\over m}}}
  \sum_{u\in (P)}{{c_{\ell}(u)c_q(u)}\over u}
 ={1\over {\varphi(q)}}\cdot \1_{\ell=q}, 
$$
\par
\noindent
following from the \lq \lq Smooth Twisted Orthogonality\rq \rq, see Proposition 2 in 3rd version of [C1] : 
$$
{1\over {\varphi(\ell) {\displaystyle \sum_{m\in (P)} }{1\over m}}}
  \sum_{u\in (P)}{{c_{\ell}(u)c_q(u)}\over u}=\1_{\ell=q}, 
$$
\par
\noindent
whence by Lemma 7 the thesis.\hfill \square 
\medskip
\par
\noindent {\bf Remark 13.} Once fixed $q'\in \N$, in case $\exists \chi\neq \chi_0(\bmod q')$, writing $\asymp$ for both $\ll$ and $\gg$, 
as $P\in \Primes$, $P\to \infty$, 
$$
\sum_{m\in (P)}{1\over m}=\prod_{p\le P}\left(1-{1\over p}\right)^{-1} \asymp \log P,
\enspace \hbox{\rm while } 
\sum_{m\in (P)}{{\chi(m)}\over m}=\prod_{p\le P}\left(1-{{\chi(p)}\over p}\right)^{-1} \asymp 1,
\forall \chi\neq \chi_0(\bmod q'), 
$$
\par
\noindent
as these last products converge, for $P\in \Primes$, $P\to \infty$, to \enspace $\prod_{p}(1-\chi(p)/p)^{-1}=L(1,\chi)\neq 0$. These two partial Euler products, from Lemma 9, will appear in next Theorem, in the way its sketchy Proof suggests.\hfill $\diamond$

\medskip
\par
\noindent {\bf Remark 14.} By the way, more precisely, next Theorem's bound comes from (now, $\forall q'\in \N$ fixed): 
$$
d|n \Rightarrow \left|\prod_{p|d}\left(1-\overline{\chi}(p)\right)\right|=\prod_{p|d}\left|1-\overline{\chi}(p)\right|
     \le 2^{\omega(d)}
      \le 2^{\omega(n)}
       \le 2\varphi(n), 
\quad \forall \chi(\bmod \, q'), 
$$
\par
\noindent
because : $\forall n\in \N$, 
$$
{{2^{\omega(n)}}\over {\varphi(n)}}=\prod_{p|n}{2\over {\varphi(p^{v_p(n)})}}
=\prod_{2|n}{2\over {2^{v_p(n)-1}}}\cdot \prod_{{p|n}\atop {p>2}}{2\over {(p-1)p^{v_p(n)-1}}}\le 2, 
$$
\par
\noindent
an absolute constant.\hfill $\diamond$ 
\vfill
\eject

\par 				
\noindent
We are ready to state and prove our most interesting result about Carmichael coefficients, both smooth and classic, for the imaginary exponential function. We may abbreviate $P-$Carmichael Transform to mean: $P-$smooth Carmichael Transform. Also, \lq \lq to\rq \rq, hereafter, may shorten \lq \lq converges to\rq \rq. 
\smallskip
\par
\noindent {\bf Theorem 5}. ({\stampatello Imaginary exponentials' $P-$Carmichael Transform to Carmichael Transform})
\par
\noindent
{\it Fix $q\in \N$ and $j\in \Z_q^*$, choose $P\in \Primes$ with $P\ge q$ and take $\ell \in (P)$. Then, the} {\stampatello explicit formula holds:} 
$$
\Carmichael_{\ell}^{(P)}\, e_q(j\bullet)=\Carmichael_{\ell}\, e_q(j\bullet) 
+ \sum_{b|q}{1\over {b\varphi(q')}}\sum_{{\chi\neq \chi_0}\atop {(\!\!\bmod q')}}\tau\left(\overline{\chi}\right)\chi(j)
  {{\mu(q'')\chi(\ell'')}\over {\varphi(\ell')}} \prod_{p|\ell''}\left(1-\overline{\chi}(p)\right)
   {{ {\displaystyle \prod_{p\le P} } \left(1-{{\chi(p)}\over p}\right)^{-1}}\over { {\displaystyle \prod_{p\le P} } \left(1-{1\over p}\right)^{-1}}}\, , 
$$ 
\par
\noindent
{\stampatello abbreviating } $q':=q/b$, $\ell':={{\ell}\over {(\ell,b)}}$, $q'':=\prod_{p|\ell',p|q'}p^{v_p(\ell')}$ {\it and } $\ell'':=\ell'/q''$. {\it As a consequence, } {\stampatello the bound:} 
$$
\Carmichael_{\ell}^{(P)}\, e_q(j\bullet)=\Carmichael_{\ell}\, e_q(j\bullet)+O_q\left({1\over {\log P}}\right)
 =\1_{\ell=q}\cdot {1\over {\varphi(q)}}+O_q\left({1\over {\log P}}\right),
$$
\par
\noindent
{\stampatello uniformly } $\forall \ell \in \N$ ({\it see Remark 14}), {\stampatello with the constant depending at most on the fixed } $q\in \N$.
\smallskip
\par
\noindent {\bf Proof(Sketch)}.Gather: Lemmas 9,10, Kluyver Formula for $c_{\ell'}(t)$, Lemma 11,12 and Remarks 13,14.\hfill \square 
\medskip
\par
Since any Correlation, say $C_{f,g_Q}(N,a)$, satisfying Basic Hypothesis is a linear combination of imaginary exponentials $e_q(ja)$ as follows: 
$$
C_{f,g_Q}(N,a)=\sum_{q\le Q}\widehat{g_Q}(q)\sum_{j\in \Z_q^*}S_f\left({j\over q}\right)e_q(ja), 
$$
\par
\noindent
where we'll abbreviate henceforth
$$
S_f(\alpha)\defineq \sum_{n\le N}f(n)e(n\alpha),
\quad \forall \alpha \in [0,1], 
$$
\par
\noindent
previous Theorem for imaginary exponentials has the following Corollary for $\BH-$correlations, of two fixed $f,g_Q:\N \rightarrow \C$. (For the details about truncated $g_Q$ and its Ramanujan coefficients $\widehat{g_Q}$, see the above $\S1.2$.)  
\smallskip
\par
\noindent {\bf Corollary 4}. ({\stampatello All $\BH-$correlations' $P-$Carmichael Transform to Carmichael Transform})
\par
\noindent
{\it Fix $Q,N\in \N$, with $Q\le N$, and abbreviate $F(a):=C_{f,g_Q}(N,a)$, $\forall a\in \N$, for the $f$ and $g_Q$ \hfil $\BH-$correlation. Choose $P\in \Primes$ and take $\ell \in (P)$. Then, the} {\stampatello explicit formula holds:} 
$$
\Carmichael_{\ell}^{(P)}\, F=\Carmichael_{\ell}\, F 
+\sum_{q\le Q}\widehat{g_Q}(q)\sum_{b|q}{1\over {b\varphi(q')}}\sum_{{\chi\neq \chi_0}\atop {(\!\!\bmod q')}}\tau\left(\overline{\chi}\right)
  \sum_{j\in \Z_q^*}\chi(j)S_f\left({j\over q}\right)\times
$$
$$
\times
   {{\mu(q'')\chi(\ell'')}\over {\varphi(\ell')}} \prod_{p|\ell''}\left(1-\overline{\chi}(p)\right)
    {{ {\displaystyle \prod_{p\le P} } \left(1-{{\chi(p)}\over p}\right)^{-1}}\over { {\displaystyle \prod_{p\le P} } \left(1-{1\over p}\right)^{-1}}}\, , 
$$ 
\par
\noindent
{\stampatello abbreviating } $q':=q/b$, $\ell':={{\ell}\over {(\ell,b)}}$, $q'':=\prod_{p|\ell',p|q'}p^{v_p(\ell')}$ {\it and } $\ell'':=\ell'/q''$. {\it As a consequence, } {\stampatello the bound:} 
$$
\Carmichael_{\ell}^{(P)}\, F=\Carmichael_{\ell}\, F+O_{Q,N,f,g}\left({1\over {\log P}}\right)
 =\widehat{g_Q}(\ell)\sum_{j\in \Z_{\ell}^*}S_f\left({j\over {\ell}}\right){1\over {\varphi(\ell)}}+O_{Q,N,f,g}\left({1\over {\log P}}\right),
$$
\par
\noindent
{\stampatello uniformly } $\forall \ell \in \N$, {\stampatello with an absolute constant depending at most on the fixed } $Q,N\in \N,$ $f,g\in {\C}^{\N}$.
\medskip
\par
\noindent {\bf Remark 15.} Since $\sum_{j\in \Z_{\ell}^*}S_f\left({j\over {\ell}}\right)=\sum_{n\le N}f(n)c_{\ell}(n)$, the RHS explicit part here is nothing but the $\ell-$th coefficient in REEF's RHS. We explicitly highlight : CONVERGENCE OF COEFFICIENTS DOESN'T IMPLY CONVERGENCE OF EXPANSIONS! (Compare Counterexample 1 in $\S5.6$, for this.)\hfill $\diamond$ 

\vfill
\eject

\par 				
\noindent
Of course we might consider a kind of {\stampatello Approximate Reef} for $\BH-$correlations, once defined the {\stampatello Error Term}: 
$$
E_{f,g_Q}(N,a)\defineq C_{f,g_Q}(N,a) - \sum_{q\le Q}\widehat{g_Q}(q)\sum_{n\le N}f(n)c_q(n){1\over {\varphi(q)}}c_q(a) 
$$
\par
\noindent
arising (as we saw in previous Corollary and Remark) as our $\BH-$correlation minus its REEF RHS. With this definition, previous Corollary may be written more explicitly, for the Correlation, as
\smallskip
\par
\noindent {\bf Corollary 5}. ({\stampatello Explicit Formula For Error Term of $\BH-$correlations})
\par
\noindent
{\it Let $F(a):=C_{f,g_Q}(N,a)$, $\forall a\in \N$, represent a $\BH-$correlation and define its Error Term as above. Then $\forall a\in \N$ fixed, choosing $P\in \Primes$ with $P\ge \max(Q,a)$, we get:} 
$$
E_{f,g_Q}(N,a)=\sum_{\ell\in (P)}c_{\ell}(a)
   \left(\sum_{q\le Q}\widehat{g_Q}(q)\sum_{b|q}{1\over {b\varphi(q')}}\sum_{{\chi\neq \chi_0}\atop {(\!\!\bmod q')}}\tau\left(\overline{\chi}\right)
    \sum_{j\in \Z_q^*}\chi(j)S_f\left({j\over q}\right)
\times
\right.
$$
$$
\left.
\times
   {{\mu(q'')\chi(\ell'')}\over {\varphi(\ell')}} \prod_{p|\ell''}\left(1-\overline{\chi}(p)\right)
    {{ {\displaystyle \prod_{p\le P} } \left(1-{{\chi(p)}\over p}\right)^{-1}}\over { {\displaystyle \prod_{p\le P} } \left(1-{1\over p}\right)^{-1}}}\right), 
$$ 
\par
\noindent
{\it where this quantity in brackets is $\Carmichael_{\ell}^{(P)}\,F-\Carmichael_{\ell}\,F=\Wintner_{\ell}^{(P)}\,F-\Wintner_{\ell}\,F$ and we } {\stampatello abbreviate as above } $q':=q/b$, $\ell':={{\ell}\over {(\ell,b)}}$, $q'':=\prod_{p|\ell',p|q'}p^{v_p(\ell')}$ {\it and } $\ell'':=\ell'/q''$. 
\medskip
\par
We leave the Proof as an exercise, for the interested reader. 
\medskip
\par
As we also leave the other following \lq \lq Exercise\rq \rq, arising from the question: what if we introduce Dirichlet characters AT ONCE from the imaginary exponential $e_q(ja)$ in the $\BH-$correlation? 
\smallskip
\par
\noindent {\bf Theorem 6}. ({\stampatello Dirichlet Characters Explicit Formula For $\BH-$correlations Error Term})
\par
\noindent
{\it Let $F(a):=C_{f,g_Q}(N,a)$, $\forall a\in \N$, represent a $\BH-$correlation and define its Error Term as above. Then $\forall a\in \N$ fixed, abbreviate now $q':=q/(q,a)$ and $a':=a/(q,a)$, to get:} 
$$
E_{f,g_Q}(N,a)=\sum_{q\le Q}\widehat{g_Q}(q)\cdot {1\over {\varphi(q')}}\sum_{{\chi\neq \chi_0}\atop {(\!\!\bmod q')}}\tau\left(\overline{\chi}\right)
    \chi(a')\sum_{j\in \Z_q^*}\chi(j)S_f\left({j\over q}\right). 
$$

\vfill
\eject

\par 				
\noindent{\bf 6. Odds \& ends. Recent work. Further remarks \& future work}
\bigskip 
\par
\noindent
We start with some complementary results, about Euler products, in next subsection $\S6.1$.
\par
Then, subsection $\S6.2$ of remarks on the connections between $F'$ and $\WintnerT F$, for general $F:\N \to \C$. In present Version 9 it is expanded, from new results described in next subsection. 
\par
Recent work, leading to a few considerable new results, is in fact given in $\S6.3$.
\par
\noindent
Last but not least, we give further remarks and a kind of \lq \lq short coming soon\rq \rq \enspace on future work: last subsection $\S6.4$. 

\bigskip
\bigskip

\par
\noindent{\bf 6.1. Euler products}
\bigskip 
\par
\noindent
We give a very short proof of a property coming from the equivalence in Theorem 2. 
\medskip
\par
We have its immediate application to Euler products, even if we need further hypotheses. 
\smallskip
\par
\noindent {\bf Proposition 4}. {\it If $F:\N \rightarrow \C$ has $\WintnerT F$, is multiplicative, satisfying}
$$
F'=\mu^2\cdot F',
\quad 
\lim_x \sum_{{r\in )P(}\atop {1<r\le x}}{{F'(r)}\over r}=\lim_x \sum_{{r\in )P(\cap (x)}\atop {r>1}}{{F'(r)}\over r}, 
\enspace \forall P\in \Primes, 
\enspace \hbox{\it and} 
\quad 
\sum_p \log\left(1+{{F'(p)}\over p}\right)\enspace \hbox{\it converges}, 
$$ 
\par
\noindent
{\it then $F$ has an Eratosthenes transform with finite support: $|\supporto(F')|<\infty$.}
\smallskip
\par
\noindent {\bf Proof(Sketch)}. From the equivalence of Theorem 2, using the hypothesis over the limits for $x\to \infty$, we can express the $r-$series as an infinite Euler product : (the case $F=\0$ has a trivial proof, so we know $F\neq \0$, that implies $F'(1)=F(1)=1$ here) 
$$
\sum_{r\in )P(}{{F'(r)}\over r}=\prod_{p>P}\Big(1+{{F'(p)}\over p}\Big) 
 =\exp\Big(\sum_{p>P}\log\Big(1+{{F'(p)}\over p}\Big)\Big) 
  \buildrel{P}\over{\to} 1, 
$$
\par
\noindent
from the hypothesis of convergence for the $\log-$series over primes.\qed 
\medskip
\par
An immediate application to $F(a):=C_{f,g}(N,a)$, thanks to $\BH$ consequences (see above), gives our 
\smallskip
\par
\noindent {\bf Corollary 6}. {\it If $C_{f,g}(N,a)$ satisfies $\BH$ and the following  hypotheses: $C_{f,g}(N,\cdot)$ is multiplicative, with $C_{f,g}'(N,\cdot)$ square-free supported},
$$
\lim_x \sum_{{r\in )P(}\atop {1<r\le x}}{{C_{f,g}'(N,r)}\over r}=\lim_x \sum_{{r\in )P(\cap (x)}\atop {r>1}}{{C_{f,g}'(N,r)}\over r}, 
\enspace \forall P\in \Primes, 
\quad 
\hbox{\it and}
\quad
\sum_p \log\left(1+{{C_{f,g}'(N,p)}\over p}\right)\enspace \hbox{\it converges}, 
$$ 
\par
\noindent
{\it then $C_{f,g}(N,a)$ has the } $\REEF$. 
\medskip
\par
\noindent
Correlations like in $\HL$, of course, are not multiplicative. However, future work can be devoted to this specific hypothesis, for general $F$ (namely, applying Proposition 4). 

\bigskip
\bigskip

\par
\noindent{\bf 6.2. Eratosthenes Transforms and their averages (following Wintner)}
\bigskip
\par
\noindent
We give here some remarks, about $F'$ and $\WintnerT F$ links. 
\bigskip
\par
\noindent
By definition, from Eratosthenes transform (always existing), when \enspace $\exists \WintnerT F$, we know Wintner transform. 
\par
Even in case \enspace $\exists \WintnerT F$, on the other hand, we can't identify $F'$ from the knowledge of \enspace $\WintnerT F$. 
\par
\noindent
Actually, this is not completely true: with delicate assumptions, our formula \enspace $(7)$ \enspace can help, say, \lq \lq to rebuild $F'$ from $\WintnerT F$\rq \rq. Namely, the Wintner Orthogonal Decomposition for $F'$ helps knowing $F'$ from $\WintnerT F$. 
\par
Wintner Orthogonal Decomposition comes from a kind of \lq \lq arithmetic orthogonality\rq \rq: 
$$
\forall P\in \Primes, \forall d\in \N,
\quad
d=d_{(P)}\cdot d_{)P(}, 
\quad
\hbox{where} \enspace d_{(P)}=\prod_{p\le P}p^{v_p(d)} \enspace \thinspace \hbox{and} \thinspace \enspace d_{)P(}\defineq \prod_{p>P}p^{v_p(d)} 
$$
\par				
\noindent
are, say, the $P-$smooth, resp., the $P-$sifted part of $d$ (and $d_{(P)}$ defined in $(5)$ above, with the usual $p-$adic valuation recalled soon before Lemma 1). 
\par
From this, in fact, once we consider (compare Theorem 1 proof in $\S1$) 
$$
\sum_{d\not \in (P)}{{F'(d)}\over d}\sum_{{q\in (P)}\atop {q|d}}c_q(a) = \sum_{{d\not \in (P)}\atop {d_{(P)}|a}}{{F'(d)}\over d}\cdot d_{(P)} 
 = \sum_{{d_{)P(}>1}\atop {d_{(P)}|a}}{{F'(d_{(P)}\cdot d_{)P(})}\over {d_{)P(}}}
  = \sum_{d|a} \sum_{{r\in )P(}\atop {r>1}}{{F'(dr)}\over r}, 
$$
\par
\noindent
it is clear that we are separating $P-$smooth indices, involving Wintner transform, from $P-$sifted indices, involving Eratosthenes transform. 
\par
Once we know Wintner coefficients, philosophically speaking (say, without assumptions), in order to know $F'$ at a fixed $d\in \N$, we only require knowledge of our $F'$ at natural numbers with \lq \lq arbitrarily large\rq \rq \thinspace prime factors. 
\par
In fact, compare $\S5.4$, the Wintner Assumption for $F$, thanks to Theorem 1, allows to calculate not only $F$, but also $F'$. In some sense, $\WA$ constraint on $\WintnerT F$ allows to \lq \lq rebuild\rq \rq, say, $F'$ from $\WintnerT F$. However we can not do this, in general, since two functions with the same Wintner Transform may differ a lot: for example, both $\0$ and the error term for $\BH-$correlations, see $5.7$, have Wintner Transform $\0$, but these error terms are not always $\0$ (the null-function), as testified above in $\S5.6$. 

\bigskip

In this present version 9, we add new, recent work.

\bigskip

\par
\noindent
We have found another requirement on $F$ that allows to recover $F'$ from $\WintnerT \,F$. 

\medskip

\par
It is not as easy as $\WA$ above, since it involves TWO HYPOTHESES, on our $F$ : the FIRST regards, so to speak, the SMOOTHNESS OF $\WintnerT \,F$, while the SECOND is, little by little, more and more technical (from next Theorem 7 to Theorem 8 and Theorem 9) and may be called a kind of VERTICAL CONSTRAINT so to speak. In fact, see the following, we start asking (for 2nd hypothesis) an easy condition: $F$ $\IPP$, see Th.m 7; then, after an easy definition before Th.m 8, we ask more generally that $F'$ is supported over numbers $d$ with prime-power factors $p^j$ having $j\le K$ for a {\it fixed} $K$ (generalizing previous condition $K=1$), $K$ in natural numbers, and we express this saying that $F'$ has VERTICAL LIMIT $K\in \N$, see Th.m 8; then, we ask an even MORE GENERAL condition on $F$, while keeping $\WintnerT \,F$ smooth-supported, in Theorem 9: a kind of VERTICAL CONSTRAINT that involves the Irregular Series of our $F$. 

\medskip

We explicitly warn the reader that we give Theorems 7,8,9 in order of increasing generality, to keep a kind of \lq \lq historic discovery order\rq \rq, so to speak. Also, our exposition starts from easier second hypothesis, keeping first hypothesis constant, for a kind of clarity unfolding, as concepts become more and more general. The final Theorem 9 being most general, it has Theorem 8 as a Corollary, whereas Theorem 7 is a kind of particular case ($K=1$) of Theorem 8 (general $K\in \N$), then. 
\par
Finally, see that, actually the Vertical Constraint, Theorem 9 second hypothesis, is rather cumbersome and I guess not so easy to check, see the comments soon after Theorem 9 Proof. 

\bigskip

\par
\noindent
We give a short coming soon of next subsection results, because we wish to underline that, as above, they are in the spirit of, so to speak, rebuilding $F'$ from $\WintnerT \,F$. In fact, Theorems 7 to 9 are, actually, able to imply the $P_0-$smoothness of $F'$ support from that of $\WintnerT \,F$ support. 
\par
They do it, somehow, \lq \lq {\stampatello Crossing}\rq \rq, so to speak, the {\stampatello Horizontal Limit} on Wintner Transform, since prime-factors of moduli $q$ with $\Wintner_q\,F\neq 0$ are $p\le P_0$, {\stampatello and} the {\stampatello Vertical Limit} (not on $\WintnerT \,F$ support, but) on $F'$ support, since all prime-powers factors $p^j$ of divisors $d$ with $F'(d)\neq 0$ have $j\le K\in \N$, fixed. Informally speaking, next Theorems realize, say, a kind of \lq \lq Wintner's Crossing Property\rq \rq\thinspace ! 
 
\vfill
\eject

\par				
\noindent{\bf 6.3. Arithmetic functions' vertical limits and smooth-supported $\WintnerT\, F$ entail the $\Reef$ for $F$}
\bigskip
\par
\noindent
We start with the easiest vertical limit for $F'$ (it's square-free supported) of our $F:\N \rightarrow \C$, when $F$ $\IPP$; this, together with the hypothesis that Wintner coefficients for $F$ vanish outside the $P_0-$smooth numbers, for a certain fixed $P_0\in \Primes$, gives the $\Reef$ : see next Theorem 7. 
\par
Then, we keep this hypothesis on $\WintnerT \,F$, while generalizing the vertical limit, from $K=1$ corresponding to $F$ $\IPP$, to general $K\in \N$: see the definitions, soon after next Theorem 7 \& we apply them in its generalization, Theorem 8. Even this Theorem is actually, technically speaking, a Corollary of subsequent Theorem 9; that generalizes the concept of vertical limit, through a kind of vertical constraint, not expressed in terms of prime-powers limits, but assuming a technical convergence condition, on $F$ Irregular Series. 
\medskip
\par
Since we are going to assume the same \lq \lq {\stampatello horizontal limit}\rq \rq, say, on the Wintner coefficients in all of our subsequent results, we profit to give, in next Proposition, the resulting properties of our irregular series for $F$, that we'll use in all of next results' Proofs. 
\medskip
We recall that the following \lq \lq $P-$stability\rq \rq \enspace property has already been exposed in previous sections (esp., compare Remark 6) and follows immediately from $(7)$ in Lemma 2, like \lq \lq $P-$switching\rq \rq \enspace too. 
\smallskip
\par
\noindent {\bf Proposition 5}. ({\stampatello $P-$stability \& $P-$switching for ${\rm Irr}^{(P)}\, F$, from $\WintnerT\, F$ horizontal limit})
\par
\noindent 
{\it Let $F:\N \rightarrow \C$ have } $\WintnerT F$ {\stampatello smooth-supported}, {\it namely } 
$$
\exists P_0\in \Primes: \supporto(\WintnerT F)\subseteq (P_0). 
\leqno{\WIN_{P_0}}
$$
\par
\noindent 
{\it Then}
$$
\forall P\in \Primes, P\ge P_0, 
\enspace
{\rm Irr}^{P}_d\, F={\rm Irr}^{P_0}_d\, F, 
\qquad \hbox{\stampatello uniformly}\enspace \forall d\in \N.
$$
\par
\noindent 
{\it Furthermore, this $P-${\stampatello stability} can be combined with the other property, say, $P-${\stampatello switching}, next:}
$$
{\rm Irr}^{P_0}_d\, F=-F'(d), 
\qquad \hbox{\stampatello uniformly}\enspace \forall d\not \in (P_0).
$$
\medskip
\par 
Our first result follows, to get the $\Reef$. We avoid trivial case: $F$ constant. 
\smallskip
\par
\noindent {\bf Theorem 7}. {\it Let non-constant $F:\N \rightarrow \C$ have $\supporto(\WintnerT\, F)\subseteq (P_0)$, for some $P_0\in \Primes$, assuming $F$ $\IPP$. Then, }\enspace $\forall a\in \N$, $F(a)={\displaystyle \sum_{q\in (P_0)}\left(\Wintner_q \,F\right)c_q(a), }$ {\it whence the } $F-\Reef$.
\smallskip
\par
\noindent {\bf Proof}. We start quoting Wintner Orthogonality Decomposition for $F'$, namely $(7)$ in Lemma 2: $\forall P\in \Primes$, 
$$
F'(d)=d\sum_{s\in (P)}\mu(s)\Wintner_{ds}\, F - \IrrPdF,
\quad
\forall d\in \N, 
$$
\par
\noindent
whence \enspace $\WIN_{P_0}$ \enspace in Proposition 5 gives the \lq \lq $P-${\stampatello stability}\rq \rq \enspace of $\IrrPF$ from $P=P_0$ onwards: 
$$
\IrrPF = \Irr^{(P_0)}\, F,
\quad
\forall P\ge P_0 \enspace (P\in \Primes). 
\leqno{\rm (\ast)}
$$
\par
\noindent
Then, numbering consecutive primes from $P_0$ onwards as : $P_0<P_1<P_2<\cdots<P_m<\cdots$, 
$$
\forall d\in \N, 
\enspace 
\Irr^{(P_1)}_{d}\,F = \Irr^{(P_0)}_{d}\,F 
 =\sum_{{r\in )P_1(}\atop {r>1}}{{F'(dr)}\over r}+\sum_{{r\in )P_1(}\atop {r>1}}{{F'(dP_1r)}\over {P_1r}}
  = \Irr^{(P_1)}_{d}\, F + {1\over {P_1}}\Irr^{(P_1)}_{P_1d}\, F,
$$
\par
\noindent
(we used here $F$ $\IPP$, ignoring $P_1$ prime powers), whence, since our hypothesis \enspace $\WIN_{P_0}$ \enspace gives, from both properties in Proposition 5, \enspace $\Irr^{(P_1)}_{P_1 d}\, F = -F'(P_1 d)$, we get: $F'(P_1 d)=0$, $\forall d\in \N$; iterating on $m\in \N$, in the same way
$$
\forall m\in \N, 
\forall d\in \N, 
\enspace 
\Irr^{(P_m)}_{d}\,F = \Irr^{(P_{m-1})}_{d}\,F 
  = \Irr^{(P_m)}_{d}\, F + {1\over {P_m}}\Irr^{(P_m)}_{P_m d}\, F,
$$
\par				
\noindent
again $(\ast)$ and Proposition 5 give:
$$
\forall m\in \N, 
\forall d\in \N, 
\enspace 
F'(P_m d)= -\Irr^{(P_m)}_{P_m d}\, F =0,
\quad
\hbox{\rm entailing}
$$
$$
t\not\in (P_0)
\enspace \Rightarrow \enspace 
\exists m\in \N : t=P_m d
\quad \hbox{\rm gives} \enspace 
F'(t)=0. 
$$
\par
\noindent
In other words, we have proved that: \enspace $\supporto(F')\subseteq (P_0)$, whence 
$$
\forall a\in \N,
\quad 
F(a)=\sum_{{d\in (P_0)}\atop {d|a}}F'(d)
 =\sum_{d\in (P_0)}{{F'(d)}\over d}\sum_{q|d}c_q(a)
  =\sum_{q\in (P_0)}\left(\Wintner_q \,F\right)c_q(a),
$$
\par
\noindent
the $\Reef$ following from: $F \, \IPP$ $\Rightarrow$ $F'=\mu^2\cdot F'$ $\Rightarrow$ $\WintnerT F=\mu^2\cdot \WintnerT F$ $\Rightarrow$ $\left|\supporto(\WintnerT F)\right|\le 2^{\pi(P_0)}$.\hfill \qed 

\medskip

See that, apart from the properties that come only from the smooth support of our $\WintnerT F$, the other property of our irregular series we are applying, here, is a kind of recursion which simplifies a lot, from the other hypothesis, namely no prime-power-factors in $F'$ support!
\par
In fact, without a specific hypothesis on our $F$, this recursion is not so simple. Before we generalize $\IPP$ arithmetic functions, we give a Lemma to show how this general recursion goes, for the irregular series. By the way, we need a hypothesis for this series to converge, namely: existence of Wintner Transform.

\medskip

\par
This is, so to speak, contained in next Lemma. The Proof comes from $\Irr^{(P)}\,F$ definition. 
\smallskip
\par 
\noindent {\bf Lemma 13}. ({\stampatello recursion for the irregular series})
\par
\noindent
{\it Let } $F:\N \rightarrow \C$ {\it have } $\WintnerT \,F$. {\it Then, given a sequence of consecutive primes } $P_0<P_1<\cdots<P_m<\cdots$, {\it once fixed any $m\in \N$}, 
$$
\forall d\in \N,
\quad
\Irr^{(P_{m-1})}_d \,F = \Irr^{(P_m)}_d \,F + \sum_{{r\in )P_m(}\atop {r>1}}\sum_{j=1}^{\infty}{{F'(dP^{j}_m r)}\over {P^{j}_m r}}. 
$$

\par
\noindent
Notice : the $r-$series and the $j-$series may not be exchanged, in general. Furthermore, any bound on $F'$ modulus in this double series ruins the convergence of present $r-$series ! However, if the $j-$summation is finite, we can exchange summations very easily : for this reason, we introduce a generalization of $\IPP$ functions, that have $j\le 1$, to $j\le K$, with fixed $K\in \N$, here. 

\medskip

\par
We write, $\forall n\in \N$, 
$$
V(n)\defineq \max\{ v_p(n) : p\in \Primes\}
$$
\par
\noindent
for the, say, (global) Valuation of a natural $n\in \N$. Then, by abuse of notation, we use the same symbol for the (global) Valuation of {\bf any non-zero} arithmetic function $G:\N \rightarrow \C$, $G\neq \0$, which might be infinite this time: 
$$
V(G)\defineq \sup\{ V(n) : n\in \supporto(G) \}
$$ 
\par
\noindent
and we call $G:\N \rightarrow \C$ a $\KVL$ arithmetic function, when this $\sup$ is finite : $V(G)\in \N_0$ (the case $V(G)=0$ holding IFF the only non-zero value of $G(n)$ is at $n=1$), 
$$
G \enspace \KVL \enspace \definiz V(G)\in \N_0. 
$$
\par
\noindent
Here, $\KVL$ abbreviates \lq \lq $K-$Vertically Limited\rq \rq, as we may write for these functions : $V(G)=K$. For example, $F$ $\IPP$ if and only if : $V(F')\le 1$ (i.e., $F'$ is square-free supported): recall, $V(F')=0$ exactly for constant $F=F(1)\neq \0$. However, this $F'$ is $\KVL$ but how do we call the corresponding $F$ ? Well, it Ignores Prime Powers, being \lq \lq $>K-$th powers independent\rq \rq, recalling $K=1$ for $F$ $\IPP$, and we introduce, say, the \lq \lq $K-$Vertically Independent\rq \rq \enspace arithmetic functions 
$$
F \enspace \KVI \enspace \definiz \enspace F' \enspace \KVL. 
$$

\medskip

\par				
From previous, next Lemma. Again, for the Proof recall $\Irr^{(P)}\,F$ definition. We avoid $F$ constant, now. 
\smallskip
\par 
\noindent {\bf Lemma 14}. ({\stampatello recursion for the irregular series of $K-$vertically independent arith. fun.s})
\par
\noindent
{\it Let } $F:\N \rightarrow \C$ {\it have } $\WintnerT \,F$ {\it and let $F$ be } $\KVI$, {\it with } $V(F')=K\in \N$. {\it Then, given a sequence of consecutive primes } $P_0<P_1<\cdots<P_m<\cdots$, {\it once fixed any $m\in \N$}, 
$$
\forall d\in \N,
\quad
\Irr^{(P_{m-1})}_d \,F = \Irr^{(P_m)}_d \,F + \sum_{j\le K}{{\Irr^{(P_m)}_{dP_m^j} \,F}\over {P^{j}_m}}. 
$$

\bigskip

\par
Now, we generalize hypothesis $F$ $\IPP$ to hypothesis $F$ $\KVI$, here, with $F$ non-constant. 
\smallskip
\par
\noindent {\bf Theorem 8}. {\it Let $F:\N \rightarrow \C$ have $\supporto(\WintnerT F)\subseteq (P_0)$, for a certain $P_0\in \Primes$, and assume $F$ $\KVI$, {\it with } $V(F')=K\in \N$. Then, }\enspace $\forall a\in \N$, $F(a)={\displaystyle \sum_{q\in (P_0)}\left(\Wintner_q \,F\right)c_q(a), }$ {\it whence the } $F-\Reef$.
\smallskip
\par
\noindent {\bf Proof}. The case of second hypothesis : non-constant $F$ $\IPP$ is equivalent to $V(F')=1=K$. We start for next cases $K\ge 2$, getting the {\stampatello same consecutive primes} $P_0<P_1<\cdots<P_m<\cdots$, together with the $P-${\stampatello stability of} $\Irr^{(P)}\, F$, from $P=P_0$ on; but now, for general $K\in \N$ we need Lemma 14, to get :  
$$
\forall m\in \N, 
\forall d\in \N, 
\enspace 
\Irr^{(P_0)}_{d}\,F = \Irr^{(P_0)}_{d}\, F + \sum_{j\le K}{{\Irr^{(P_m)}_{dP_m^j} \,F}\over {P^{j}_m}}
\enspace \Rightarrow \enspace 
\sum_{j\le K}{{F'(dP_m^j)}\over {P^{j}_m}}=0
\enspace \Rightarrow \enspace 
\sum_{j\le K}{{F'(dP_m^j)}\over {P^{j-1}_m}}=0, 
$$
\par
\noindent
after using $P-$stability and $P-$switching of Proposition 5, from the hypothesis $\WIN_{P_0}$. 
\par
Thus 
$$
\forall m\in \N, 
\forall d\in \N, 
\enspace 
F'(dP_m)= -\sum_{j\le K-1}{{F'(dP_m^{j+1})}\over {P^{j}_m}},
\leqno{(\ast\ast)}
$$
\par
\noindent
after renaming the $j-$variable, here. This $(\ast\ast)$ is {\stampatello recursion on $P_m-$powers}. In fact, fix $m\in \N$ and this recursion, together with $V(F')=K$, say, \lq \lq {\stampatello kills powers}\rq \rq \enspace {\stampatello from the highest}: 
$$
F'(dP_m)= -\sum_{j\le K-1}{{F'(dP_m^{j+1})}\over {P^{j}_m}}, \forall d\in \N
\enspace (\hbox{\stampatello set}\, d:=P_m^{K-1}t) \Rightarrow \enspace 
F'(tP_m^{K})= -\sum_{j\le K-1}{{F'(tP_m^{j+K})}\over {P^{j}_m}}=0, \forall t\in \N, 
$$
\par
\noindent
whence 
$$
F'(dP_m)= -\sum_{j\le K-2}{{F'(dP_m^{j+1})}\over {P^{j}_m}}, \forall d\in \N
\enspace (d:=P_m^{K-2}t) \Rightarrow 
F'(tP_m^{K-1})= -\sum_{j\le K-2}{{F'(tP_m^{j+K-1})}\over {P^{j}_m}}=0, \forall t\in \N, 
$$
\par
\noindent
where this time we combine $V(F')=K$ with {\stampatello previous vanishing above}. Iterating, we get  
$$
F'(tP_m^2)= -{{F'(dP_m^3)}\over {P_m}}=0, \forall t\in \N
$$
\par
\noindent
from vertical limit and previous vanishing values, whence 
$$
F'(dP_m)= -{{F'(dP_m^2)}\over {P_m}}=0, \forall d\in \N. 
$$
\par
\noindent
In all, $P_m|t$ $\Rightarrow$ $F'(t)=0$ and this holds $\forall m\in \N$.
\par 
In other words, we get back to previous Proof last part, as \enspace $\supporto(F')\subseteq (P_0)$ \& so on: the $\Reef$'s from \enspace $\left|\supporto(\WintnerT F)\right|\le (K+1)^{\pi(P_0)}$.\hfill \qed 

\medskip

\par
\noindent {\bf Remark 16.} We see the irony of fate at work on $(\ast\ast)$, as {\stampatello a posteriori} it becomes completely {\stampatello  trivial}.\hfill $\diamond$ 

\medskip

\par				
Next, we give present, most general hypothesis on $F'$, here. 
\smallskip
\par
\noindent {\bf Theorem 9}. {\it Let $F:\N \rightarrow \C$ have $\supporto(\WintnerT F)\subseteq (P_0)$, for a certain $P_0\in \Primes$, and} 
$$
\forall P\ge P_0 , \forall d\in \N,
\quad
\sum_{{r\in )P(}\atop {r>1}}\sum_{j=1}^{\infty}{{F'(dP^j r)}\over {P^j r}}=\sum_{j=1}^{\infty}P^{-j}{\rm Irr}^{(P)}_{dP^j}\,F. 
$$
\par
\noindent
{\it Then, }\enspace $\supporto(F')\subseteq (P_0)$, {\it whence } {\stampatello the Ramanujan-Wintner Smooth Expansion}.
\smallskip
\par
\noindent {\bf Proof}. The Lemma 13 above gives, together with Proposition 5 like in previous Proofs, with the {\stampatello same consecutive primes} $P_0<P_1<\cdots<P_m<\cdots$, fixing $m\in \N$, after changing $j$ variable, 
$$
\forall d\in \N, 
\enspace 
\Irr^{(P_0)}_{d}\,F = \Irr^{(P_0)}_{d}\, F + \sum_{j=1}^{\infty}{{\Irr^{(P_m)}_{dP_m^j} \,F}\over {P^{j}_m}}
\enspace \Rightarrow \enspace 
F'(dP_m)=-\sum_{j=1}^{\infty}P^{-j}_m F'(dP_m^{j+1}). 
$$
\par
Thus 
$$
F'(dP_m)=-P_m^{-1}F'(dP_m^2)-\sum_{j=1}^{\infty}P^{-j-1}_m F'(dP_m^{j+2}),
\enspace \forall d\in \N, 
\leqno{(\ast\ast\ast)}
$$
\par
\noindent
whence, setting $d=tP_m$ and back with $d$ instead of $t$, 
$$
F'(dP_m^2)=-P_m^{-1}F'(dP_m^3)-\sum_{j=1}^{\infty}P_m^{-j-1} F'(dP_m^{j+3}), 
\enspace 
\forall d\in \N, 
$$
\par
\noindent
which we plug into $(\ast\ast\ast)$ to get 
$$
F'(dP_m)=-P_m^{-1}\left(-P_m^{-1}F'(dP_m^3)-\sum_{j=1}^{\infty}P_m^{-j-1} F'(dP_m^{j+3})\right)-\sum_{j=1}^{\infty}P_m^{-j-1} F'(dP_m^{j+2})
 =0, 
\enspace 
\forall d\in \N, 
$$
\par
\noindent
true $\forall m\in\N$, whence \enspace $\supporto(F')\subseteq (P_0)$.\hfill \qed 

\medskip

\par
Notice : the condition on exchanging double summation in the double series, say, in Theorem 9 \enspace \lq \lq Vertical Constraint\rq \rq, is very technical and doesn't allow easy shortcuts, as the double series doesn't converge absolutely due to the lack of absolute convergence for the Irregular series ! 

\medskip

Also, this most general result doesn't supply the $\Reef$ because it has no explicit request on $F'$ vertical LIMIT : this, in Theorems 8,7 allows to estimate explicitly the cardinality of non-vanishing Wintner coefficients, whence the $\Reef$. (Compare Theorems 8,7 Proofs final parts.) 

\par
Going back to applications, for $\BH-$correlations $F(a):=C_{f,g_Q}(N,a)$, $\forall a\in \N$, see that the condition $\WIN_{P_0}$ follows from $\BH$ (compare $\S1.2$ above), with $P_0\defineq \max\{ p\in \Primes : p\le Q\}$ and this, thanks to Counterexample 1 studied in $\S5.6$ above, renders cristal clear that we need a kind of vertical constraint. In fact, since no $\REEF$ holds for it (see quoted $\S5.6$) we see that, not only it is not $\IPP$ (compare $\S5.6$, Curiosity 1), but it has neither the much lighter vertical constraint, in Theorem 9 above. 

\vfill
\eject

\par				
\noindent
\par
\noindent{\bf 6.4. Further remarks and future work}
\bigskip
\par
\noindent
We were looking, in previous versions, for a kind of \lq \lq supplementary hypothesis\rq \rq, which, added to $\BH$, gives the $\REEF$ : we thought $\ETD$ could be the right one. Actually, $\BH$ alone doesn't give the $\REEF$, as we proved in third version of [C1] with Counterexample 1 there, compare $\S5.6$. As our Theorem 1 shows, $\WA$ is a good hypothesis of this kind : it \lq \lq gives the $\REEF$\rq \rq, to $\BH-$correlations (in Corollary 1). Finally, the \lq \lq missing hypothesis\rq \rq, say, is given by the vertical constraints, more and more general, of above Theorems 7,8,9 : in fact, from $\BH$ we know that $\supporto(\WintnerT F)\subseteq [1,Q]\subseteq (P_0)$, with $Q\le P_0\in \Primes$, abbreviating with $F(a)$ our correlation of shift $a\in \N$. Actually, Theorem 9 generality doesn't supply the $\REEF$, but only the Ramanujan-Wintner Smooth Expansion. Our Theorem 8 and its particular case Theorem 7, here, give the $\REEF$ to $\BH-$correlations, but at a high price so to speak: a vertical limit on the divisors $d|a$ of correlation's shift $a\in \N$. This is not so natural, for a correlation; however, it points in the \lq \lq heuristically right direction\rq \rq, say, i.e.: $\BH$ correlations with shift-factor $g$ which is $\IPP$ have square-free supported Wintner Transforms, entailing that $\REEF$'s main term (that's $A_F$, see $\S5.4$) is $\IPP$ itself, with smooth-supported Wintner Transform, in full concordance with Theorem 7, say!

\bigskip
\bigskip
\bigskip

\par
\noindent
There are two main directions where to look at in future work: {\bf the $\BH-$correlations world}, both for its own sake \& for the inspiration for finding (as we did in present work!) new general results; {\bf and} the theoretical fascination coming from the \lq \lq new Ramanujan clouds\rq \rq: mainly {\bf the Ramanujan smooth clouds}, as for Ramanujan clouds we already started, with Luca Ghidelli, a kind of structural description (beginning with multiplicative Ramanujan coefficients, compare [CG1] and [CG2]). 
\par
Last but not least we will, in future papers, give other explicit formul\ae, for the correlations satisfying $\BH$, coming from the elementary approach (compare Theorem 6 above) with the so-called {\bf Dirichlet} characters {\bf explicit formul\ae}. 

\vfill
\eject

\par				
\centerline{\stampatello Bibliography}

\bigskip

\item{[C0]} G. Coppola, {\sl An elementary property of correlations}, Hardy-Ramanujan J. {\bf 41} (2018), 65--76.
\smallskip
\item{[C1]} G. Coppola, {\sl A smooth shift approach for a Ramanujan expansion}, ArXiV:1901.01584v3 (Third Version)
\smallskip
\item{[C2]} G. Coppola, {\sl Finite and infinite Euler products of Ramanujan expansions}, ArXiV:1910.14640v2 (Second Version) 
\smallskip
\item{[C3]} G. Coppola, {\sl Recent results on Ramanujan expansions with applications to correlations}, Rend. Sem. Mat. Univ. Pol. Torino {\bf 78.1} (2020), 57--82. 
\smallskip
\item{[CG1]} G. Coppola and L. Ghidelli, {\sl Multiplicative Ramanujan coefficients of null-function}, ArXiV:2005.14666v2 (Second Version) 
\smallskip
\item{[CG2]} G. Coppola and L. Ghidelli, {\sl Convergence of Ramanujan expansions, I [Multiplicativity on Ramanujan clouds]}, ArXiV:1910.14640v1 
\smallskip
\item{[CM]} G. Coppola and M. Ram Murty, {\sl Finite Ramanujan expansions and shifted convolution sums of arithmetical functions, II}, J. Number Theory {\bf 185} (2018), 16--47. 
\smallskip
\item{[D]} H. Davenport, {\sl Multiplicative Number Theory}, 3rd ed., GTM 74, Springer, New York, 2000. 
\smallskip
\item{[De]} H. Delange, {\sl On Ramanujan expansions of certain arithmetical functions}, Acta Arith., 31 (1976), 259--270.
\smallskip
\item{[HL]} G.H. Hardy and J.E. Littlewood, {\sl SOME PROBLEMS OF 'PARTITIO NUMERORUM'; III: ON THE EXPRESSION OF A NUMBER AS A SUM OF PRIMES.} Acta Mathematica {\bf 44} (1923), 1--70. 
\smallskip
\item{[K]} J.C. Kluyver, {\sl Some formulae concerning the integers less than $n$ and prime to $n$}, Proceedings of the Royal Netherlands Academy of Arts and Sciences (KNAW), 9(1):408--414, 1906. 
\smallskip
\item{[M]} M. Ram Murty, {\sl Ramanujan series for arithmetical functions}, Hardy-Ramanujan J. {\bf 36} (2013), 21--33. Available online 
\smallskip
\item{[R]} S. Ramanujan, {\sl On certain trigonometrical sums and their application to the theory of numbers}, Transactions Cambr. Phil. Soc. {\bf 22} (1918), 259--276.
\smallskip
\item{[ScSp]} W. Schwarz and J. Spilker, {\sl Arithmetical Functions}, Cambridge University Press, 1994.
\smallskip
\item{[T]} G. Tenenbaum, {\sl Introduction to Analytic and Probabilistic Number Theory}, Cambridge Studies in Advanced Mathematics, {46}, Cambridge University Press, 1995. 
\smallskip
\item{[W]} A. Wintner, {\sl Eratosthenian averages}, Waverly Press, Baltimore, MD, 1943. 

\bigskip
\bigskip
\bigskip

\par
\leftline{\tt Giovanni Coppola - Universit\`{a} degli Studi di Salerno (affiliation)}
\leftline{\tt Home address : Via Partenio 12 - 83100, Avellino (AV) - ITALY}
\leftline{\tt e-mail : giocop70@gmail.com}
\leftline{\tt e-page : www.giovannicoppola.name}
\leftline{\tt e-site : www.researchgate.net}

\bye